\documentclass[11pt]{article}
\usepackage{amsmath}
\usepackage{graphicx}
\usepackage{epsfig}
\usepackage{amsfonts}
\usepackage{amssymb}
\usepackage{placeins}
\usepackage{cases}
\usepackage[margin=1in]{geometry}
\usepackage{authblk}
\usepackage[titletoc,toc,title]{appendix}
\usepackage{epstopdf}
\usepackage[final]{changes}
\usepackage{enumerate}

\usepackage{accents}

\title{An Adaptive Phase-Amplitude Reduction Framework Without $\mathcal{O}(\epsilon)$ Constraints on Inputs}

\begin{document}
\author[1]{Dan Wilson}
\affil[1]{Department of Electrical Engineering and Computer Science, University of Tennessee, Knoxville, TN 37996, USA}
\maketitle

\begin{abstract}
Phase reduction is a well-established technique used to analyze the timing of oscillations in response to weak external inputs.  In the preceding decades, a wide variety of results have been obtained for weakly perturbed oscillators that place restrictive limits on the magnitude of the inputs or on the magnitude of the time derivatives of the inputs.  By contrast, no general reduction techniques currently exist to analyze oscillatory dynamics in response to arbitrary, large magnitude inputs and comparatively very little is understood about these strongly perturbed limit cycle oscillators.  In this work, the theory of isostable reduction is leveraged to develop an adaptive phase-amplitude transformation that does not place any restrictions on the allowable input.   Additionally, provided some of the Floquet multipliers of the underlying periodic orbits are near-zero, the proposed method yields a reduction in dimension comparable to that of other phase-amplitude reduction frameworks.  Numerical illustrations show that the proposed method accurately reflects synchronization and entrainment of coupled oscillators in regimes where a variety of other phase-amplitude reductions fail.  
\end{abstract}

\maketitle


\section{Introduction}
\textcolor{black}{
Phase reduction \cite{erme10}, \cite{izhi07}, \cite{kura84} is a commonly used strategy to study oscillatory dynamical systems of the form 
\begin{equation} \label{maineq}
\dot{x} = F(x) + U(t)
\end{equation} in terms of a  1-dimensional phase reduced equation 
\begin{equation} \label{predapx}
\dot{\theta} = \omega + Z(\theta) \cdot U(t).
\end{equation}
In the full dynamical model \eqref{maineq} above, $x \in \mathbb{R}^N$ is the state, $F$ gives the nominal dynamics, and $U(t)$ is an external input; in the phase reduced equation \eqref{predapx}, $\theta \in \mathbb{S}^1$ is a phase which characterizes the timing of an oscillation, $\omega$ is the natural frequency, $Z(\theta)$ is a phase response curve, and `$\cdot$' denotes the dot product.}    Phase reduction is particularly effective when all Floquet exponents of the underlying periodic orbit are large in magnitude so that the periodic orbit behaves like an inertial manifold \cite{foia88}.  However, if any Floquet exponents are $\mathcal{O}(1)$, then some $\mathcal{O}(\epsilon)$ constraints must be placed on $U(t)$  (where $0<\epsilon \ll 1$).    For instance, phase reduction and phase-amplitude reduction \cite{wils17isored}, \cite{shir17}, \cite{cast13}, \cite{wedg13}, \cite{lets20}, \cite{bres18}, \cite{wils18operat} require the inputs to remain $\mathcal{O}(\epsilon)$ in magnitude so that the state remains close to a nominal limit cycle.  High order reduction strategies have been proposed that allow for larger magnitude inputs \cite{wils19prl}, \cite{rose19}, \cite{leon19}, \cite{wils20highacc} but these strategies still break down as the state moves far from the nominal limit cycle.  Other strategies allow for arbitrarily large magnitude inputs, but place $\mathcal{O}(\epsilon)$ limits on the rate at which the inputs can vary  \cite{kure13}, \cite{park16}, \cite{pyra15}, \added{\cite{rubi13}}.   

While phase reduction and related methods have been applied fruitfully to weakly perturbed oscillators \cite{winf01}, \cite{erme10}, \cite{piet19} there are currently no reduction frameworks that allow for the study of strongly perturbed oscillators with general perturbations.  In this work, an adaptive phase-amplitude coordinate framework is proposed that can be implemented without imposing any $\mathcal{O}(\epsilon)$  constraints on the input.  The critical feature of this reduction strategy is that it explicitly considers the relationships between phase and amplitude coordinates of nearby parameter sets.  A high-accuracy reduction results from continuously updating the underlying model parameters in order to limit the magnitude of the amplitude coordinates during the application of exogenous inputs.  Additionally, provided enough of the Floquet exponents of the periodic orbit are large in magnitude, this framework results in a significant reduction in dimensionality.    In supporting numerical examples,  this framework far exceeds the accuracy of a wide variety other reduction techniques when large magnitude inputs are considered.

The organization of this paper is as follows:~Section \ref{apxb} provides necessary background on existing phase and phase-amplitude reduced equations.  These methods are used for comparison to the adaptive phase-amplitude reduction that is the focus of this work.  Section \ref{derivsec} provides a derivation of the adaptive phase-amplitude reduced equations with Appendix \ref{apxa} giving a detailed error analysis of the truncated terms.  Section \ref{resultsec} provides illustrative examples highlighting the utility of the proposed adaptive phase-amplitude reduction when large magnitude input is considered.   In these examples, the proposed strategy faithfully reproduces entrainment and synchronization observed in simulations of the full models in situations where other strategies fail.  Section \ref{concsec} gives concluding remarks.

\section{A \added{Brief} Review of Related Phase Reduction Frameworks}  \label{apxb}
Phase transformations and phase reduction techniques are often applied to ordinary differential equations of the general form \eqref{maineq} with oscillatory dynamics.  Supposing that a model of the general form \eqref{maineq} has a stable $T$-periodic limit cycle, $x^\gamma(t)$, in many applications, it is useful to analyze the oscillatory solutions of \eqref{maineq} in terms of the timing of the oscillations.  From this perspective, taking $U(t) = 0$ one can assign a phase $\theta \in [0,2\pi)$ to all locations on the periodic orbit scaled so that $d \theta/dt = 2 \pi/T \equiv \omega$.  The definition of phase can be extended to the nonlinear basin of attraction of the limit cycle using isochrons \cite{guck75}, \cite{winf01}.  To do so, letting $\theta_1$ be the phase corresponding to $a(0) \in x^\gamma$, the $\theta_1$ level set (i.e.,~isochron) can be defined as the set of all $b(0)$ in the basin of attraction of the limit cycle such that 
\begin{equation} \label{isodef}
\lim_{t \rightarrow \infty} ||a(t) - b(t)|| = 0,
\end{equation}
where $|| \cdot ||$ denotes any vector norm.  Using the definition \eqref{isodef}, one can show that under the flow of \eqref{maineq},  $\dot{\theta} = \omega$ when $U(t) = 0$.  The phase as defined by \eqref{isodef} takes $x \mapsto \theta$,  encoding for the infinite time convergence to the periodic orbit.  Given the reduction in dimensionality, it is often useful to work in a phase reduced coordinate framework.  The following provides a \added{brief} overview of existing phase and phase-amplitude reduction frameworks that will be used as a comparison for the adaptive phase-amplitude reduction developed in this work.  \added{For a more thorough review of these and other phase reduction techniques, the interested reader is referred to both \cite{piet19} and \cite{mong19}.}

\subsection*{Standard Phase Reduction}
\added{Phase reduction \cite{erme10}, \cite{izhi07}, \cite{kura84} is a widely applied technique that can be used to represent the behavior  of a general dynamical system \eqref{maineq} near an $N$-dimensional limit cycle with the 1-dimensional representation of the form \eqref{predapx}.}  The idea of phase reduction was investigated in detail by Winfree \cite{winf01} more than half a century ago and has since become an indispensable tool for the analysis of limit cycle oscillators. A standard phase reduction can be implemented by considering the isochrons as defined by \eqref{isodef} and computing the phase response curve, i.e.,~the gradient of the isochrons evaluated on the periodic orbit.  This standard implementation of a phase reduction, as given by \eqref{predapx}, yields a particularly elegant characterization of the dynamical behavior provided the state remains close to the underlying limit cycle.  To ensure that \eqref{predapx} remains valid, it is required that $U(t)$ remain sufficiently small in magnitude relative to the decay rates of all neglected amplitude coordinates.  



\subsection*{Phase-Amplitude Reduction Methods and Isostable Coordinates}
Many authors have explored various phase-amplitude reduction strategies \cite{wils16isos},\cite{shir17},\cite{cast13},\cite{wedg13},\cite{lets20}, \cite{bres18},\cite{wils18operat},\cite{mong19b}.  Amplitude coordinates can be used in conjunction with the phase coordinates to characterize the dynamics that govern the decay in directions transverse to a periodic orbit.  This additional information can be used to better characterize the observed behavior of a perturbed oscillator \cite{wils20ddred}.  It can also be used to develop control strategies that actively limit the magnitude of the amplitude coordinates so that the state remains close to the periodic orbit.  \added{Such control strategies allow for larger magnitude inputs to be applied than those that result from using the standard phase reduction \cite{mong19b}.}  

In many applications\added{,} it can be useful to characterize the amplitude dynamics using Floquet theory \cite{jord07}.  To do so, define $\Delta x = x-x(\theta)$ to be the difference between a state $x$ with corresponding phase $\theta$ and another state $x(\theta)$ on the periodic orbit.  To a linear approximation, \eqref{maineq} can be written as
\begin{equation}
\Delta \dot{x} = J \Delta x,
\end{equation}
where $J$ is the (time-varying) Jacobian evaluated at $x^\gamma(\theta)$.  Letting $\Phi$ be the fundamental matrix yielding the relationship $\Delta x(T) = \Phi \Delta x(0)$, supposing that $\Phi$ is diagonalizable, the solutions near the periodic orbit can be represented according to \cite{wils17isored}
\begin{equation} \label{outputerror}
x - x^\gamma(\theta) = \sum_{j = 1}^{N-1} \psi_j g^j(\theta) + \mathcal{O}(\psi_1^2) + \dots + \mathcal{O}(\psi_{N-1}^2),
\end{equation}
where $g^j(\theta)$ are Floquet eigenfunctions associated with the periodic orbit $x^\gamma(\theta)$, and $\psi_j$ are isostable coordinates.  These isostable coordinates can be defined as level sets of particular eigenfunctions of the Koopman operator \cite{maur13}, \cite{kval19}.  For the most slowly decaying isostable coordinates, an alternative definition can be given that considers the infinite-time decay of initial conditions to the limit cycle \cite{wils17isored}, \cite{shir17}.  \added{To linear order,} the dynamics of the phase-isostable reduction can be written as 
\begin{align} \label{fordred}
\dot{\theta} &= \omega + Z(\theta) \cdot U(t), \nonumber \\
\dot{\psi}_j &= \kappa_j \psi_j + I_j(\theta) \cdot U(t),  \nonumber \\
&j = 1,\dots,N-1,
\end{align}
where $\kappa_j$ is a Floquet exponent associated with the $j^{\rm th}$ isostable coordinate and $ I_j(\theta)$ gives the gradient of $\psi_j$ evaluated on the periodic orbit.  Typically, $U(t)$ is assumed to be an order $\epsilon$ term so that \eqref{fordred} can be referred to as a `first order' accurate phase-amplitude reduction where all of the phase and isostable coordinate dynamics are accurate to $\mathcal{O}(\epsilon)$.  Additionally, if some of the Floquet exponents are large in magnitude so that they die out rapidly, the associated isostable coordinates will remain close to zero and can be ignored.  This leaves $\beta < N-1$ isostable coordinates to consider in the reduction \eqref{fordred}.  


\subsection*{Second and Higher Order Accurate Phase-Amplitude Reduction Methods}
While \eqref{predapx} and \eqref{fordred} only use information about the gradients of the phase and amplitude coordinates evaluated on the periodic orbit (i.e.,~computations of the gradients that are valid to zeroth order in the amplitude coordinates) computation of these gradients at higher order accuracy can illuminate model behaviors that are not replicated by  the standard phase reduced dynamics \cite{wils19prl}, \cite{rose19}, \cite{leon19}, \cite{cast13}.   A second order accurate version of the phase-isostable reduced equations was developed in \cite{wils17isored} and \cite{wils19complex}:
\begin{align} \label{firstord}
\dot{\theta} &= \omega + \left(Z(\theta)  + \sum_{k = 1}^{\beta} \psi_k B^k(\theta)    \right) \cdot U(t),     \nonumber \\
\dot{\psi}_j &= \kappa_j \psi_j +  \left(I_j(\theta)  + \sum_{k = 1}^{\beta} \psi_k C_j^k(\theta)    \right) \cdot U(t),   \nonumber \\
& j = 1, \dots, \beta.
\end{align}
The functions $B^k(\theta)$ and $C_j^k(\theta)$ in the above equation provide higher order corrections to the gradient of the phase and isostable dynamics.   Provided $U(t)$ is an $\mathcal{O}(\epsilon)$ term, Equation \eqref{firstord} is valid to second order accuracy in $\epsilon$.  Recent work \cite{wils20highacc} developed techniques to compute reduced order equations to arbitrarily high order accuracy.

\subsection*{Extended Phase Reduction}
The notion of an extended phase reduction proposed in \cite{kure13} (cf.~\cite{park16}) allows for arbitrary magnitude inputs as long as they vary sufficiently slowly.  This strategy splits the input $U(t)$ into two components: $U(t) = q(\epsilon t) + \sigma r(t)$.  Here, $q(\epsilon t)$ is a slowly varying component of arbitrary magnitude and $\sigma r(t)$ is an additional weak component with $0 <\sigma \ll 1$.   The imposed structure of $U(t)$ allows the amplitude coordinates to be ignored in the reduction yielding an extended phase reduced equation of the form \cite{kure13}
\begin{equation} \label{extpred}
\dot{\theta} = \omega(q(\epsilon t)) + D(\theta,q(\epsilon t)) \cdot \dot{q}(\epsilon t) + \sigma Z(\theta,q(\epsilon t)) \cdot r(t),
\end{equation}
where $D(\theta,q(\epsilon t))$ characterizes the response to the slowly varying parameter and $Z(\theta,\xi)$ is the phase response curve of the periodic \added{orbit that results} when a constant $U(t) = \xi$ is applied.  Equation \eqref{extpred} is valid provided that \added{the magnitudes} of the Floquet exponents of the periodic orbit are small relative to both $\dot{q}(\epsilon t)$ and $\sigma r(t)$.   The extended phase reduction \eqref{extpred} allows  for large magnitude inputs to be used provided these inputs vary slowly enough--a constraint that precludes the use of \eqref{extpred} in many applications.  
\textcolor{black}{
\subsection*{Related Reduced Order Coordinate Frameworks for Nonlinear Dynamical Systems}
  Recent years have seen a sustained interest in developing tools and techniques that can be used to characterize dynamical behavior far from an underlying attractor.  To this end, the Koopman operator paradigm has gained traction in recent years and can be used to represent a general nonlinear dynamical system using a linear, but infinite-dimensional operator \cite{budi12}, \cite{mezi13}.  The fundamental challenge of implementing Koopman-based approaches is in the identification of a finite basis with which to represent the potentially infinite dimensional Koopman operator.  As mentioned earlier, the isostable coordinate framework leverages Koopman-based ideas by working in a finite basis of the slowest decaying Koopman eigenmodes \cite{maur13}, \cite{kval19}.  Other strategies such as Dynamic  Mode Decomposition (DMD) \cite{schm10}, \cite{kutz16}, extended DMD \cite{will15}, and delayed embedding approaches \cite{arba17}, \cite{brun17} are well suited for identifying reduced order Koopman-based models, especially in data-driven applications.}

\added{
Alternatively, the parameterization method \cite{haro16} can be applied to identify invariant manifolds of a dynamical system.  This general strategy has been used as a foundation to compute sets of phase-amplitude coordinates (such as isochrons and isostables) in the basin of attraction of a limit cycle \cite{cast13}, \cite{guil09}.  Recently, numerical techniques based on the parameterization method have been developed to efficiently compute full sets of isochron and isostable coordinates as well as their gradients for dynamical systems of arbitrary dimension \cite{albe20}.}




\section{Derivation of an Adaptive Phase-Amplitude Reduction Strategy} \label{derivsec}
A derivation of the \added{proposed} adaptive phase-amplitude reduction is presented here.  In the previous section, each of the phase and phase-amplitude reduction strategies fail when the state $x$ is driven too far from the underlying periodic orbit $x^\gamma(\theta)$.  \added{The extended phase reduction \eqref{extpred} discussed in the previous section takes steps to remedy this by considering a family of periodic orbits that emerge as a given parameter is changed. However,  it does not explicitly consider the associated amplitude dynamics.  Consequently, limits must be placed on the rate of change of the magnitude of the allowable inputs. } The critical feature of the proposed adaptive phase-amplitude reduction strategy is that it actively updates the system parameters (and consequently the nominal periodic orbit) with the explicit goal of keeping the state close to the nominal periodic orbit.  

\textcolor{black}{
To begin, consider a general ordinary differential equation of the form
\begin{equation} \label{maineqadapt}
\dot{x} = F(x,p_0) + U(t).
\end{equation}
\added{where $x$, $F$, and $U$ are identical to those from Equation \eqref{maineq}, and $p_0 \in \mathbb{R}^M$ is a constant, nominal collection of system parameters. Note here that Equation \eqref{maineqadapt} is in the same form as \eqref{maineq} with the explicit inclusion of a nominal parameter set.   A shadow system}
\begin{equation} \label{shadoweq}
\dot{x} = F(x,p) + U(t),
\end{equation}  
is introduced, where $p \in \mathbb{R}^M$ is a set of adaptive parameters.  The ultimate goal is to identify a reduced order of system of equations that accurately reflect the behavior of \eqref{maineqadapt} with the help of the shadow system \eqref{shadoweq}.   It will be assumed that for all allowable $p$ that  \eqref{maineqadapt} has a periodic orbit denoted by $x_{p}^\gamma(t)$.}  \added{The notion of phase can be extended to the entire basin of attraction of each periodic orbit using isochrons as defined in Equation \eqref{isodef}.  For a given parameter set, initial conditions with the same asymptotic approach to the limit cycle are defined to have the same phase.  To this end, let $\theta(x,p)$ denote the phase associated with the limit cycle $x_{p}^\gamma(t)$ and the state $x$.}   \added{Of course, the isochrons of each periodic orbit defined according to \eqref{isodef} are unique up to a constant shift, a consideration that will be considered more carefully momentarily}.

To implement the adaptive phase-amplitude reduction framework, suppose that \eqref{maineqadapt} has a nominal parameter set $p_0$ that describes the dynamics of an underlying model.  It is then possible to rewrite \eqref{maineqadapt} as
\begin{equation} \label{extendedframe}
\dot{x} = F(x,p) + U_e(t, p, x),
\end{equation}
with the extended input  
\begin{equation} \label{extinput}
U_e(t,p,x) \equiv U(t) + F(x,p_0)-F(x,p).
\end{equation}
Intuitively, the term $F(x,p)$ represents dynamics of the model when using parameter set $p$, and $U_e$ captures both the externally applied input and the mismatch caused when considering the parameter set $p$ as compared to the parameter set $p_0$ of the underlying model \eqref{maineqadapt}.    In the following adaptive phase-amplitude reduction formulation, $p$ will be viewed as a free variable.  Changing to phase and isostable coordinates via the chain rule yields
\begin{align} \label{fulleq}
\frac{d}{dt} \theta(x,p) &= \frac{\partial \theta}{\partial x} \cdot  \frac{d x}{dt} + \frac{\partial \theta}{\partial p} \cdot   \frac{d p}{dt}, \nonumber \\
\frac{d}{dt } \psi_j(x,p) &= \frac{\partial \psi_j}{\partial x} \cdot  \frac{d x}{dt} + \frac{\partial \psi_j }{\partial p} \cdot  \frac{d p}{dt}, \nonumber  \\
j &= 1, \dots, \beta \nonumber \\
\frac{dp}{dt} &= G_p(p,\theta,\psi_1,\dots,\psi_\beta)
\end{align}
 where the function $G_p$ dictates how $p$ is updated.  \added{The above reduction truncates $\beta + 1,\dots,N-1$ rapidly decaying isostable coordinates.  Conditions describing when this truncation is possible are discussed in Section \ref{assumpsec}}.

\added{Recall that the isochrons of each $x_p^\gamma(t)$ are unique up to a constant shift; in order for the terms of the form $\frac{\partial \theta}{\partial p}$ and $\frac{\partial \psi_j}{\partial p}$ \eqref{fulleq} to be well-defined, these constant shifts must be unambiguously determined for the continuous family of periodic orbits.  In addressing a similar problem, the authors of \cite{wang19} disambiguate the phase shift between limit cycles by defining a Poincar\'e section transverse to a set of limit cycles and mandating that the phase at the intersection of this Poincar\'e section and each limit cycle be identical. An identical construction can be employed here in the definition of $\theta(x,p)$.  Likewise, one can also require a distinct feature associated with the family of limit cycles to have identical phase between limit cycles, for instance, by setting the maximum transmembrane voltage during an action potential to correspond to an identical phase on each limit cycle.  While there are multiple options, one of these strategies must be used to unambiguously define $\theta(x,p)$. }

 \added{Equation \eqref{fulleq} provides the basis for the proposed adaptive reduction framework in which changes to either the system state $x$ or the the adaptive parameter set $p$ can influence both the phase and isostable coordinates.  Most phase reduction techniques, including those highlighted in Section \ref{apxb}, rely on an asymptotic expansion of the phase in a basis of amplitude coordinates about a periodic attractor.  As such, they fail when the amplitude coordinates (e.g.,~isostable coordinates) become so too large.  Intuitively,  the adaptive parameter set $p$ provides an additional degree of freedom that can be used to keep the amplitude coordinates small when large magnitude inputs are applied.  The sections to follow formalize this intuition and provide strategies for computing each of the terms of \eqref{fulleq}.}
 
 \textcolor{black}{
\subsection{Necessary Assumptions and Restrictions on the Proposed Adaptive Reduction Framework} \label{assumpsec}
Explicit assumptions and restrictions for implementation of the adaptive reduction are given below.  These assumptions ultimately allow for simplification of the terms of \eqref{fulleq} and subsequent implementation of the proposed adaptive reduction framework.   The resulting reduction explicitly considers $\beta$ slowly decaying isostable coordinates and truncates the remaining rapidly decaying isostable coordinates.
\begin{enumerate}[{Assumption} A)]
\item In the allowable range of $p$,  the periodic orbit $x_p^\gamma(t)$ exists and is continuously differentiable with respect to both $t$ and $p$.
\item In the allowable range of $p$, $\theta(x,p)$ and $\psi_j(x,p)$ are continuously differentiable with respect to $x$ and $p$ for all $j$.
\item The function $G_p$ that governs the adaptive parameter is designed  so that for all time,  $|\psi_j|$ is an $\mathcal{O}(\sqrt \epsilon)$ term for $j \leq \beta$ where $0 < \epsilon \ll 1$.
\item For the truncated isostable coordinates, $\min_{p,j > \beta} (\kappa_j(p)) = \mathcal{O}(1/\epsilon)$ where $0 < \epsilon \ll 1$
\item The norms of $U_e(t,p,x)$ and $G_p(p,\theta,\psi_1,\dots,\psi_\beta)$ are bounded for all time and allowable $p$ so that $||U_e(t,p,x)||_1 \leq M_U$ and  $|| G_p(p,\theta,\psi_1,\dots,\psi_\beta) ||_1 \leq M_p$ where $|| \cdot ||_1$ denotes the $1$-norm.
\end{enumerate}
 Assumptions A and B ensure that the dynamical system under consideration is sufficiently smooth so that the necessary partial derivatives from \eqref{fulleq} exist.  These assumptions generally exclude piecewise smooth dynamical systems such as those considered in \cite{park16ar} and \cite{wils19complex}.  Additionally, assumptions A and B generally exclude critical points of bifurcations from the allowable parameter set;  it is emphasized that oscillations that emerge as a result of a bifurcation can still be readily considered provided the critical point is excluded from the allowable parameter set.   Assumptions C, D, and E  ensure that the resulting reduced order model remains accurate to leading order $\epsilon$ even for very large inputs.  A brief discussion about design heuristics for $G_p$ is given in Section \ref{designgp}.  A detailed error analysis of the resulting reduced order equations is presented in Appendix \ref{apxa}.
 }
 
 \textcolor{black}{
As explained in Appendix A, as a consequence of Assumptions D and E, each $\psi_j$ is an order $\epsilon$ term for $j = \beta + 1, \dots, N-1$.  Assumption C mandates that each $\psi_j$ is an order $\sqrt{\epsilon}$ term for $j = 1,\dots,\beta$.  This coupled with a first order expansion in the Floquet eigenfunctions described by \eqref{outputerror} yields
\begin{equation} \label{floqeq}
x - x_p^\gamma(\theta) = \sum_{k = 1}^{\beta}  \psi_k g^k(\theta,p) + \mathcal{O}(\psi_1^2) + \dots \mathcal{O}(\psi_{N-1}^2) = \mathcal{O}(\sqrt{\epsilon}),
\end{equation}
where $g^k(\theta,p)$ is the Floquet eigenfunction associated with the periodic orbit $x^\gamma_p$.  In other words, Assumptions C, D, and E imply that $x - x_p^\gamma(\theta_0)$ is an $\mathcal{O}(\sqrt{\epsilon})$ term.
 }

 \subsection{Computation of the Necessary Terms of the Adaptive Reduction}
\added{Below, it is illustrated how each of the partial derivatives of  \eqref{fulleq} can be written in a basis of phase and isostable coordinates.   First considering the terms, $\frac{\partial \theta}{\partial x} \cdot  \frac{d x}{dt}$ and $\frac{\partial \psi_j}{\partial x} \cdot  \frac{d x}{dt}$, in \cite{wils19complex} and \cite{wils19phase} it was shown that models of the form \eqref{extendedframe} can be represented using phase and isostable coordinates according to }
\begin{align} \label{firstordadapt}
\dot{\theta} &= \frac{\partial \theta}{\partial x} \cdot  \frac{d x}{dt} = \omega(p) + \left(Z(\theta,p)  + \sum_{k = 1}^{\beta} \psi_k B^k(\theta,p)    \right) \cdot U_e(t),     \nonumber \\
\dot{\psi}_j &= \frac{\partial \psi_j}{\partial x} \cdot  \frac{d x}{dt} = \kappa_j(p) \psi_j +  \left(I_j(\theta,p)  + \sum_{k = 1}^{\beta} \psi_k C_j^k(\theta,p)    \right) \cdot U_e(t),   \nonumber \\
& j = 1, \dots, \beta.
\end{align}
Equation \eqref{firstordadapt} is valid to first order accuracy \added{in the non-truncated isostable coordinates}; the truncated isostable coordinates decay rapidly so that they can be neglected. Consequently, if $U(t)$ is an $\mathcal{O}(\epsilon)$ term, then \eqref{firstordadapt} is valid to second order accuracy in $\epsilon$.   Once again, \eqref{firstordadapt} is identical to \eqref{firstord} except that the terms of the reduced order equations depend explicitly on $p$.  Note here that $\dot{\theta}$ and $\dot{\psi}_j$ from \eqref{firstordadapt} are identical to $\frac{d \theta}{dt}$ and $\frac{d \psi_j}{dt}$, respectively, when $p$ is held constant.  The terms of right hand sides of \eqref{firstordadapt} can be computed for a chosen value of $p$ using methods described in \cite{wils19complex}. 

 
   \textcolor{black}{
For the remaining terms from \eqref{fulleq} of the form $\frac{\partial \theta}{\partial p}$ and $\frac{\partial \psi_j}{\partial p}$, consider an arbitrary initial condition $x_0$ for which $\theta(x_0,p) = \theta_0$.   Recall from \eqref{floqeq} that $x_0-x^\gamma_p$ is an $\mathcal{O}(\sqrt{\epsilon})$ term.  Taking an asymptotic expansion centered at the periodic orbit, for any state $x$ in the neighborhood $x_0$ the phase and each isostable coordinates to are given by
\begin{align} \label{thetaequation}
 \theta(x,p) &= \theta_0 + (x-x_p^\gamma) \cdot Z(\theta_0,p) +  \frac{1}{2} (x-x_p^\gamma)^T H_\theta(\theta_0,p)  (x-x_p^\gamma) + \mathcal{O} (\epsilon^{3/2}), \\ 
 \psi_j(x,p) &=  (x-x_p^\gamma) \cdot I_j(\theta_0,p) +  \frac{1}{2} (x-x_p^\gamma)^T H_{\psi_j}(\theta_0,p)  (x-x_p^\gamma) + \mathcal{O} (\epsilon^{3/2}),  \label{psiequation}
 \end{align}
for $j = 1,\dots,\beta$, where $H_{\theta}(\theta_0,p)$ (resp.,~$H_{\psi_j}(\theta_0,p)$) is the Hessian matrix of second partial derivatives of $\theta$ (resp.,~$\psi_j$) evaluated at $\theta_0$ and $p$ on the limit cycle, and $x_p^\gamma$ is evaluated at $\theta_0$.  Taking the partial of \eqref{thetaequation} with respect to $p_i$ yields
\begin{align} \label{firstpartial}
\frac{\partial \theta}{\partial p_i} &=-\frac{\partial x_p^\gamma}{\partial p_i} \cdot Z(\theta_0,p) + \frac{\partial Z}{\partial p_i} \cdot(x-x_p^\gamma)  \nonumber \\
&\quad - \frac{1}{2} \left(  \frac{\partial x_p^\gamma}{\partial p_i}^T H_\theta(\theta_0,p)  (x-x_p^\gamma)  -  (x-x_p^\gamma)^T \frac{H_\theta}{\partial p_i} (x-x_p^\gamma)  +  (x-x_p^\gamma)^T H_\theta(\theta_0,p)  \frac{\partial x_p^\gamma}{\partial p_i }\right),
\end{align}
where $^T$ denotes the vector transpose and all partial derivatives are evaluated at $\theta_0$ and $p$.  Noting that the Hessian is symmetric, neglecting order $\epsilon$ terms and simplifying yields
\begin{equation} \label{heqn}
\frac{\partial \theta}{\partial p_i} =  -\frac{\partial x_p^\gamma}{\partial p_i} \cdot \left( Z(\theta_0,p) + H_\theta(\theta_0,p)  (x-x_p^\gamma) \right) +  \frac{\partial Z}{\partial p_i} \cdot(x-x_p^\gamma) + \mathcal{O}(\epsilon).
\end{equation}
Finally,  substituting Equation \eqref{floqeq} into \eqref{heqn} one finds
\begin{align} \label{heqn}
\frac{\partial \theta}{\partial p_i} &=  -\frac{\partial x_p^\gamma}{\partial p_i} \cdot \left( Z(\theta_0,p) + \sum_{k = 1}^{\beta}  \psi_k  H_\theta(\theta_0,p)   g^k(\theta_0,p) \right) +   \sum_{k = 1}^{\beta}  \psi_k  \frac{\partial Z}{\partial p_i} \cdot g^k(\theta_0,p) + \mathcal{O}(\epsilon) \nonumber \\
&=  -\frac{\partial x_p^\gamma}{\partial p_i} \cdot \left( Z(\theta_0,p) + \sum_{k = 1}^{\beta}  \psi_k B^k(\theta,p) \right)  +   \sum_{k = 1}^{\beta}  \psi_k  \frac{\partial Z}{\partial p_i} \cdot g^k(\theta_0,p) + \mathcal{O}(\epsilon),
\end{align}
where the relationship $B^k(\theta,p) = H_\theta(\theta_0,p)  g^k(\theta_0,p) $ used in the second line was illustrated in \cite{wils19complex}.  Evaluating  \eqref{heqn} for each adaptive parameter and collecting terms appropriately, one can write
\begin{align} \label{dtdps}
\frac{\partial \theta}{\partial p} &= D(\theta,p) + \sum_{k = 1}^{\beta} \psi_k E^k(\theta,p)  +  \mathcal{O}(\epsilon),
 \end{align}
 where $D(\theta,p)$ and $E^k(\theta,p)$ are periodic in $\theta$.  In an analogous manner, starting by taking partial derivatives of \eqref{psiequation} with respect to $p_i$ and mirroring the steps that yield equations \eqref{firstpartial}--\eqref{heqn}, one finds
 \begin{equation} \label{hpsieqn}
\frac{\partial \psi_j}{\partial p_i}  =  -\frac{\partial x_p^\gamma}{\partial p_i} \cdot \left( I_j(\theta_0,p) + \sum_{k = 1}^{\beta}  \psi_k C_j^k(\theta,p) \right)  +   \sum_{k = 1}^{\beta}  \psi_k  \frac{\partial I_j}{\partial p_i} \cdot g^k(\theta_0,p) + \mathcal{O}(\epsilon).
\end{equation}
 Once again, evaluating \eqref{hpsieqn} for each adaptive parameter and collecting terms appropriately, one can write,
  \begin{align} \label{dpdps}
  \frac{\partial \psi_j}{\partial p} =  Q_j(\theta,p) + \sum_{k = 1}^{\beta} \psi_k R_j^k(\theta,p) + \mathcal{O}(\epsilon),
 \end{align}
where $Q_j(\theta,p)$ and $R_j^k(\theta,p)$ are periodic in $\theta$.  
 }
Both \eqref{dtdps} and \eqref{dpdps} are valid to first order accuracy in the non-truncated isostable coordinates.    Substituting \eqref{dtdps}, \eqref{dpdps}, and \eqref{firstordadapt} into \eqref{fulleq} yields the adaptive phase-amplitude transformation
\begin{align} \label{truncmain}
\dot{\theta} &= \omega(p) + \left(Z(\theta,p)  + \sum_{k = 1}^{\beta} \psi_k B^k(\theta,p)    \right) \cdot U_e(t,p,x)   \nonumber \\
 &+ \left( D(\theta,p) + \sum_{k=1}^{\beta} \psi_k E^k(\theta,p)    \right) \cdot \dot{p} + \mathcal{O}(\epsilon),   \nonumber \\
\dot{\psi}_j &= \kappa_j(p) \psi_j +  \left(I_j(\theta,p)  + \sum_{k = 1}^{\beta} \psi_k C_j^k(\theta,p)    \right) \cdot U_e(t,p,x)   \nonumber \\
&+  \left( Q_j(\theta,p) + \sum_{k=1}^{\beta} \psi_k R_j^k(\theta,p)     \right) \cdot \dot{p} + \mathcal{O}(\epsilon), \nonumber \\
j &= 1, \dots, \beta, \nonumber \\
\dot{p} &= G_p(p,\theta,\psi_1,\dots, \psi_\beta),
\end{align}
where $p$ evolves according to the function $G_p$. Notice that only order $\epsilon$ terms are truncated from the adaptive reduction \eqref{truncmain}; the error analysis for the truncated terms is presented in detail in Appendix \ref{apxa}.  For a given value of $p$, Equation \eqref{floqeq} provides an estimate of the state to leading order in the isostable coordinates.   Considering the aforementioned order of accuracy of the isostable coordinates, the relationship
\begin{equation} \label{firstorderdx}
 x  = x^\gamma_p(\theta) +  \sum_{j = 1}^{\beta} \psi_j g^j(\theta,p)
\end{equation}
provides an order $\epsilon$ accurate approximation for the state.  


 This reduction \eqref{truncmain} actively sets the nominal parameters $p$ so that the system state is always close to the periodic orbit $x^\gamma_p$.   In order to implement this reduction, the exact choice of $G_p$ that determines $\dot{p}$ is not important as long as it keeps the non-truncated isostable coordinates small.  \added{While the addition of the $p$ dynamics in \eqref{truncmain} adds to the overall dimension of the reduction, there are no $\mathcal{O}(\epsilon)$ restrictions on the input $U(t)$ (for instance in its magnitude or the magnitude of the first derivatives). }
 

\textcolor{black}{
\subsection{General Heuristics For Designing $\mathbf{G_p}$}\label{designgp}
Here, the design of $G_p$ in order to limit the magnitude of the non-truncated isostable coordinates (as required by Assumption C given in Section \ref{assumpsec}) is discussed.   First, consider a situation where both $p \in \mathbb{R}$ and only one non-truncated isostable coordinate, $\psi_1$, is required.  In this case, consider a function of the form 
\begin{equation} \label{firstgp}
G_p(p, \theta, \psi_1) = -\alpha \psi_1  \left( Q_1(\theta,p) +  \psi_1 R_1^1(\theta,p) \right),
\end{equation}
 where $\alpha > 0$.  Substituting Equation \eqref{firstgp} into $\frac{d \psi_1}{d t}$ from Equation \eqref{truncmain} and simplifying yields
\begin{align} 
\frac{d \psi_1}{d t} &=  \kappa_1(p) \psi_1 +  \left(I_j(\theta,p)  +  \psi_1 C_1^1(\theta,p)    \right) \cdot U_e(t,p,x)    -\alpha \psi_1  \left( Q_1(\theta,p) +  \psi_1 R_1^1(\theta,p) \right)^2 \nonumber \\
&= \left( \kappa_1(p)  -\alpha  ( Q_1(\theta,p) +  \psi_1 R_1^1(\theta,p) )^2 \right)\psi_1 +  \left(I_j(\theta,p)  +  \psi_1 C_1^1(\theta,p)    \right) \cdot U_e(t,p,x).
\end{align}
Above, the term $-\alpha  ( Q_1(\theta,p) +  \psi_1 R_1^1(\theta,p) )^2 < 0$ and serves to drive the magnitude of the isostable coordinate to lower values.  This general form of $G_p$ is used in the examples from Sections \ref{nonradsec} and \ref{neursec}.  A similar formulation of $G_p$ is applied in Section \ref{circsec} where it is adapted to a situation with 3 non-truncated isostable coordinates.}
 
 \textcolor{black}{
As an alternative formulation, once again consider an application where $p \in \mathbb{R}$ and only one non-truncated isostable coordinate, $\psi_1$ is required.  Supposing $|Q_1(\theta,p) +  \psi_1 R_1^1(\theta,p)|$ is sufficiently bounded away from zero for all $\theta, p$, and $\psi_1$, one can take
\begin{equation}\label{canceleqn}
G_p(p,\theta,\psi_1) = \frac{ \left( I_j(\theta,p)  +  \psi_1 C_1^1(\theta,p) \right) \cdot U_e(t,p,x) }{Q_1(\theta,p) +  \psi_1 R_1^1(\theta,p)}.
\end{equation}
Substituting Equation \eqref{canceleqn} into $\frac{d \psi_1}{d t}$ from Equation \eqref{truncmain} and simplifying yields
\begin{equation}  \label{psisimplified}
\frac{d \psi_1}{dt} = \kappa_1(p) \psi_1.
\end{equation}
Since \eqref{psisimplified} has a globally stable equilibrium at $\psi_1 = 0$, the $\psi_1$ dynamics can be ignored from the adaptive reduction formulation.  Such an approach was implemented in \cite{wils21adapt} to investigate optimal control inputs for influencing oscillation timing.  }

\added{Equations \eqref{firstgp} and \eqref{canceleqn} are meant to provide intuition about the design of $G_p$ but are not straightforwardly generalizable  to applications where more than one adaptive parameter is used or where more than one non-truncated isostable coordinate is required.    The general design of the parameter update rule $G_p(p,\theta,\psi_1,\dots,\psi_\beta)$ could be posed as a nonlinear control problem and will be the subject of future investigations.}

\textcolor{black}{
\subsection{Choosing the Adaptive Parameter Set}} \label{adaptchoose}
\added{The ability to design an adequate $G_p$ (as discussed in Section \ref{designgp}) is directly related to the choice of the adaptive parameter set $p$ used in the adaptive reduction framework.  Ultimately, the implementation of the adaptive reduction requires $x - x_p^\gamma(\theta)$ to be an order $\sqrt{\epsilon}$ term.  As such $p$ must be chosen so the set of $x_p^\gamma(\theta)$ is sufficient accomplish this task.}

\added{As a general heuristic for the choice of $p$, it often works to choose the adaptive parameter set such that it mirrors the effect of an external input $U(t)$.  For instance, in the example from Section \ref{neursec} that considers the reduction of a capacitance based-neuron, the external inputs are capable of influencing the model through direct transmembrane current injections.  Subsequently, the adaptive parameter is taken to be a constant transmembrane current that mirrors the influence of these external inputs.  This general strategy is also applied in the examples from Sections \ref{nonradsec} and \ref{circsec}.  }

\added{As a final note, in each of the examples considered in the work, the inputs are of the form $U(t) = \eta u(t)$ where $\eta \in \mathbb{R}^N$ is fixed and $u(t) \in \mathbb{R}$.  In other words, the inputs considered are of low rank relative to the dimension of the full order dynamical system.  Consequently, only a single adaptive parameter (that mirrors the influence of these rank-1 inputs) is required to implement the adaptive reduction strategy in these examples.  In general, higher rank inputs will require more than one adaptive parameter.}


\section{Illustrative Examples} \label{resultsec}

Here, the adaptive phase-amplitude reduction approach from \eqref{truncmain} is applied to a variety of numerical models.  Results are compared to the reduction techniques described in Section \ref{apxb}.  In each application, the proposed reduction strategy significantly outperforms the previously developed reduction techniques when large magnitude inputs are considered.

\subsection{Nonradial Isochron Clock} \label{nonradsec}
\begin{figure}[htb]
\begin{center}
\includegraphics[height=2.5in]{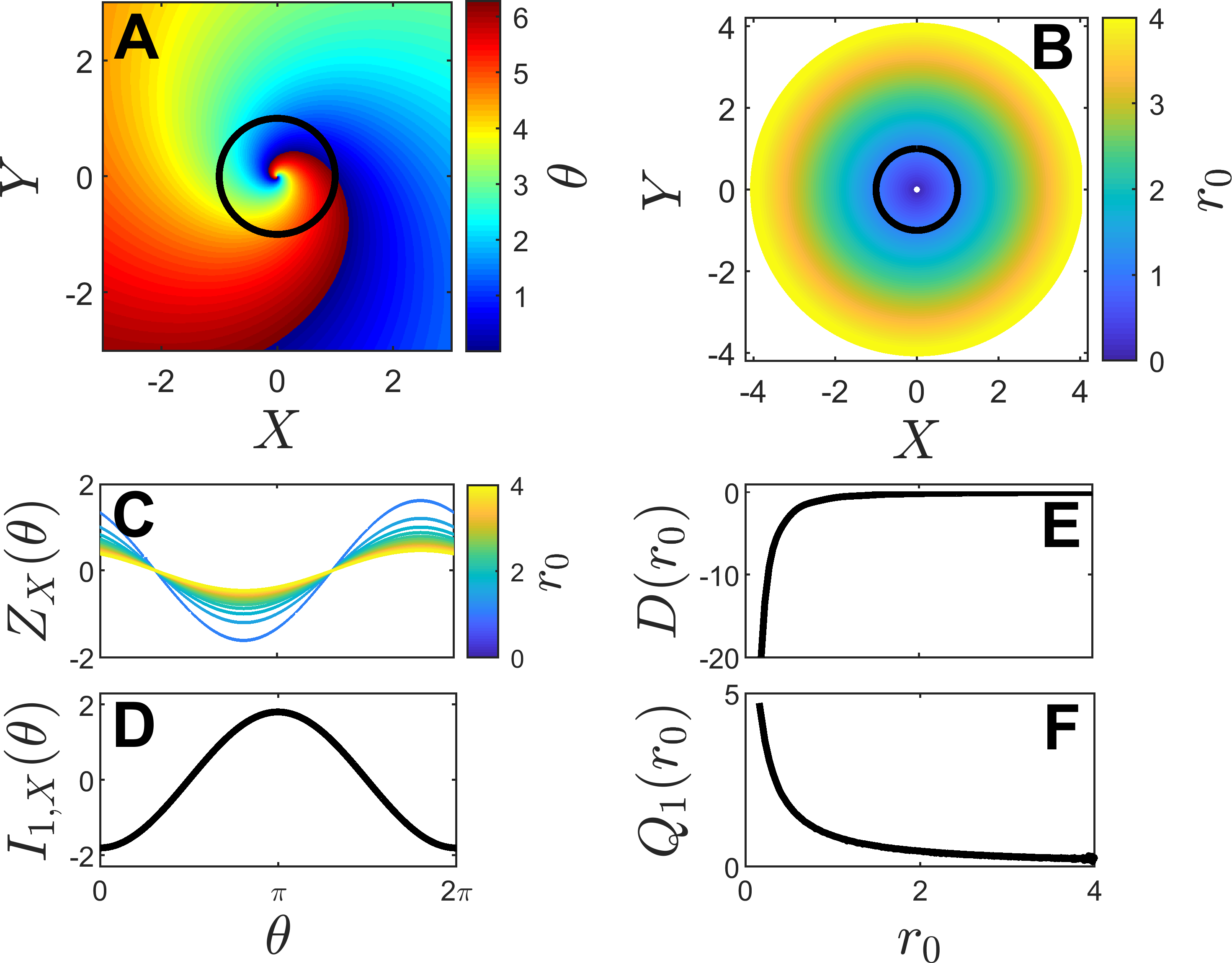}
\end{center}
\caption{Isochrons of the nonradial isochron clock are shown in Panel A.  The black line show the nominal limit cycle that emerges when $r_0 = 1$.  The radius of the limit cycle changes as $r_0$ varies according to the colormap shown in panel B.  The adaptive phase-amplitude reduction actively sets $r_0$ to keep the isostable coordinate (i.e.,the distance from the periodic orbit) small.  Panel C shows the $X$ component of the phase response curve as $r_0$ varies.  The $X$ component of the isostable response curve does not depend on $r_0$ and is shown in panel D.  Panels E and F show $D(r_0)$ and $Q_1(r_0)$ from \eqref{truncmain}.  In general these functions also depend on $\theta$, but due to radial symmetry of \eqref{isoclock} they only depend on the \added{adaptive} parameter $r_0$ in this example. }
\label{radialinfo}
\end{figure}

As a first example,  the nonradial isochron clock is considered: 
\begin{align} \label{isoclock}
\dot{X} &= \mu X (r^2_0- (X^2+Y^2)) - Y (1+\zeta( (X^2+Y^2) -r^2_0)) + u(t), \nonumber \\
\dot{Y} &= \mu Y (r^2_0- (X^2+Y^2)) + X (1+\zeta( (X^2+Y^2) -r^2_0)).
\end{align}
Above, $X$ and $Y$ are Cartesian coordinates.  For all parameter sets, the limit cycle is a circle with radius $r_0$ centered at the origin, $\mu = 0.08$ sets the convergence rate, and $\zeta = 0.12$ influences the rate of rotation. Here, \eqref{isoclock} has been modified from the radial isochron clock equations from \cite{winf01} so that the rate of rotation depends on the radius.   The parameter $r_0$ will be the \added{adaptive} parameter used in the adaptive phase-amplitude reduction, and the nominal value will be set to 1 in simulations of \eqref{isoclock}.    For the nominal parameter set, panel A of Figure \ref{radialinfo} shows the periodic orbit as a black line with corresponding isochrons calculated according to \eqref{isodef}.  These isochrons are highly nonlinear. Consequently, the standard phase reduction will only be effective for relatively small inputs so that the state remains near the limit cycle.  

The adaptive phase-amplitude reduced equations \eqref{truncmain} account for the nonlinearity in \eqref{isoclock} by continuously choosing $r_0$ to shift the nominal periodic orbit so that the isostable coordinates remain small.  The colormap in panel $B$ shows how the periodic orbit changes as $r_0$ is modified and panel C gives the the component of the phase response curve for perturbations to the parameter $X$ as $r_0$ changes.  \added{In this example, $\theta(x,p) = 0$ corresponds to the moment that both $Y = 0$ and $X>0$ on each of the limit cycles parameterized by the choice of $r_0$}. Relevant terms of the adaptive phase-amplitude reduction are shown in panels D-F.   The isostable response curve for this system does not change with $r_0$ and is shown as a function of $\theta$ in panel D.  Due to the radial symmetry of \eqref{isoclock},  $D(\theta,r_0)$ and $Q_1(\theta,r_0)$ from \eqref{truncmain} (i.e.,~the functions that describe how the isostable and phase coordinates change as $r_0$ varies) do not depend on $\theta$; the dependence on $r_0$ of these functions is shown in Panels E and F, respectively.  \added{As $r_0$ decreases towards 0, the dynamics of \eqref{isoclock} approach the critical point of a Hopf bifurcation.  As mentioned earlier, critical points of bifurcations cannot be considered in the adaptive parameter set.  Correspondingly, the magnitudes of $D(r_0)$ and $Q_1(r_0)$ begin to approach infinity as $r_0$ approaches 0.}

\begin{figure}[htb]
\begin{center}
\includegraphics[height=2.2in]{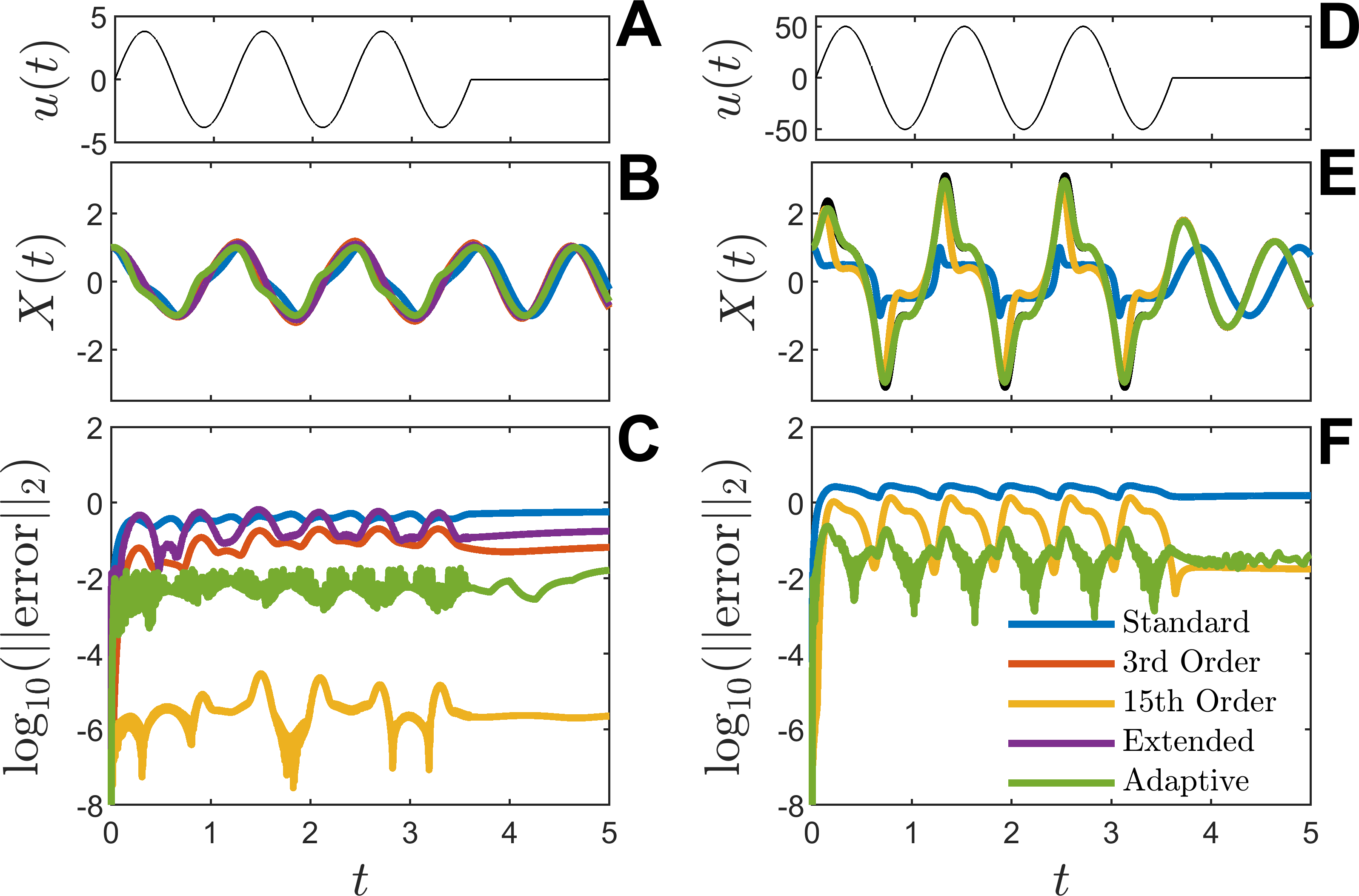}
\end{center}
\caption{Panels A and D shows small and large magnitude input, respectively, applied to the nonradial isochron clock \eqref{isoclock}.  The output from each reduced order model considered is shown in color in panels B (small magnitude input) and E (large magnitude input).  The black line shows simulations from the full model output, although it is mostly indistinguishable from the adaptive phase-amplitude reduced model.  Corresponding errors between the full and reduced model simulations are shown in panels C and F.  In panel F, some reduced order models cannot be used or do not provide viable outputs for reasons described in the text.  Unlike all other reduced models, the accuracy of the proposed adaptive phase-amplitude reduced model is not significantly degraded when using a large magnitude input.}
\label{radisomultireduc}
\end{figure}

The \added{adaptive} parameter for this model is $r_0$ with a nominal value of $r_0 = 1$, consequently the extended input from \eqref{extendedframe} is
\begin{equation}
 U_e(t,r_0,X,Y) = \begin{bmatrix}u(t) +  (\mu X + Y \zeta) -(\mu X r_0^2 + Y \zeta r_0^2)  \\  (\mu Y  - X \zeta) - (\mu Y r_0^2 - X \zeta r_0^2) \end{bmatrix}.
 \end{equation}
The adaptive phase-amplitude reduction has only one \added{adaptive} parameter, $r_0$, and one isostable coordinate, $\psi_1$.  Therefore, a relatively simple $G_p$ as suggested by \eqref{firstgp} can be chosen:~$\dot{r}_0 = G_p(r_0,\theta,\psi_1) = -\alpha \psi_1 \frac{\partial \psi_1}{\partial r_0}$ where $\alpha > 0$ and $\partial \psi_1/\partial r_0 = Q_1(\theta,r_0) + R_1^1(\theta,r_0)$.    Intuitively, this choice of $G_p$ always acts to modify the nominal parameter $r_0$ in order to decrease the isostable coordinate and updates $r_0$ more aggressively when the isostable coordinates are larger in magnitude.

The performance of the adaptive phase-amplitude reduced equations is compared to other reduction techniques discussed in Section \ref{apxb} with results shown in Figure \ref{radisomultireduc}.  Panel A shows a sinusoidal input $u(t) = 3.8 \sin(2\pi t/1.2)$ applied for 3.6 seconds.  The initial condition is taken to be $X = 1$ and $Y = 0$ (i.e.,~on the periodic orbit).  Panel B compares results when using various reduction techniques.  The standard phase reduction is given by \eqref{predapx}.  Third and fifteenth order reductions are computed using methods described in \cite{wils20highacc} by taking the gradient of the phase and isostable response curves and the reduced model output $x(\theta,\psi_1)$ to the appropriate order in $\psi_1$.    The extended phase reduction is implemented according to \eqref{extpred}.  For this reduction technique, $u(t)$ itself is considered to be a slowly varying input yielding the reduced equations  $\dot{\theta} = \omega(u(t)) + D(\theta,u(t))  \dot{u}(t)$.  For the input from panel A, each reduction strategy performs comparably well, tracking the output of the full model simulations \eqref{isoclock}.  Panel C shows the error between the full and reduced model simulations according to ${\rm error}(t) = x(t)-x^*(t)$ where $x(t) = \begin{bmatrix} X(t) & Y(t) \end{bmatrix}$ is found from the full model solutions and $x^*(t)$ is obtained from the specified reduced order model.  Note here that while the 15th order accurate model outperforms the adaptive phase-amplitude reduced model when using the logarithmic scale, they both perform reasonably well in absolute terms.  Analogous results are shown in panels D-F of Figure \ref{radisomultireduc} using $u(t) = 50 \sin(2\pi t/1.2)$.  Note that for inputs this large, the extended reduction framework \eqref{extpred} cannot be implemented, because the reduced order system does not have a stable periodic orbit for constant inputs with magnitude larger than $|u(t)| \approx 4$. Additionally, the magnitude of the input is too large for the 3rd order accurate reduction and the model output tends towards infinite values.    For the other reduction strategies that are viable when using this large magnitude input, the  adaptive phase-amplitude reduction performs significantly better.  Also, the accuracy of the adaptive phase-amplitude reduction is not significantly degraded when applying the larger magnitude input.  

It should be noted that when applied to the nonradial isochron clock \eqref{isoclock} the adaptive phase-amplitude reduction strategy yields a model with more  variables (three) than the underlying system (two).  The main goal of this illustrative example is to highlight the fact that this general framework works well even when exceedingly large magnitude inputs are used.  In the examples to follow, the adaptive phase-amplitude reduction will yield a reduced order set of equations while still providing a faithful reproduction of the unreduced model dynamics.

\subsection{Phase Locking in Synaptically Coupled Neurons} \label{neursec}
For this example, phase locking in a conductance-based model of \added{$N_{\rm neur}$} synaptically coupled neurons \cite{rubi04} will be considered:
\begin{align} \label{neurmod}
C \dot{V}_i &= -I_L(V_i)-I_{Na}(V_i,h_i)-I_K(V_i,h_i)-I_T(V_i,r_i) + I_b^i - g_{\rm syn} \sum _{j = 1}^{N_{\rm neur}} I_{\rm syn}(V_i,w_j) + u(t), \nonumber \\
\dot{h}_i &= (h_\infty(V_i)-h_i)/\tau_h(V_i), \nonumber \\
\dot{r}_i &= (r_\infty(V_i) - r_i)/\tau_r(V_i), \nonumber \\
\dot{w}_i&= \alpha(1-w_i)/(1+\exp(-(V_i-V_T)/\sigma_T)) - \beta w_i,
\end{align}
for $i = 1,\dots,N_{\rm neur}$.  Above, $V_i$ gives the transmembrane voltage of the $i^{\rm th}$ neuron, $h_i$ and $r_i$ are gating variables, and $w_i$ is a synaptic variable.   Leak, sodium, potassium, and low-threshold calcium currents, are $I_L = g_L(V_i - E_L)$, $I_{Na} = g_{Na}m_\infty^3(V_i) h_i (V_i-E_{Na})$, $I_K = g_K(.75(1-h_i))^4(V_i-E_K)$, $I_T = g_T p_\infty^2(V_i) r_i(V_i-E_T)$, respectively.  Conductances are taken to be $g_L = 0.05$, $g_{Na} = 3$, $g_K = 5$, and $g_T = 5$ ${\rm mS}/{\rm cm}^2$, reversal potentials are $E_L = -70$, $E_{Na} = 50$, $E_K = -90$, and $E_T = 0$ mV, and $C = 1 \mu {\rm F}/{\rm cm}^2$.  Parameters $V_{\rm syn} = 0$ mV, $\alpha = 3 {\rm m s}^{-1}$, $V_T = -20$ mV, $\sigma_T = 0.8$ mV, and $\beta = 2 {\rm ms}^{-1}$ set the dynamics of the synaptic variable with a synaptic coupling current given by $I_{\rm syn}(V_i,w_j) = w_j(V_i-V_{\rm syn})$.  The parameter $g_{\rm syn}$ sets the coupling strength, $I_b^i$ is a constant baseline current, \added{and $u(t)$ is a common transmembrane current input felt by each neuron}.   \textcolor{black}{Using the proposed adaptive reduction framework, each individual neuron will be reduced from a 4-dimensional model given by Equation \eqref{neurmod} to a 3-dimensional model. In the context of the reduction framework, the external input felt by each neuron is taken to be $U^i(t) = \begin{bmatrix}  (u(t) - g_{\rm syn} \sum _{j = 1}^{N_{\rm neur}} I_{\rm syn}(V_i,w_j))  &0& 0& 0  \end{bmatrix}^T $}.  For each neuron from \eqref{neurmod}, a baseline current  $I_b^i \in [2,20] \; \mu {\rm A}/\mu {\rm F}$ is treated as an \added{adaptive} parameter in the adaptive phase-amplitude reduction so that 
\begin{equation}
\textcolor{black}{U_e^i = \begin{bmatrix} u(t)  - \sum _{j = 1}^N g_{\rm syn} I_{\rm syn}(V_i,w_j) + I_{b,0}^i- I_b^i \\ 0 \\ 0 \\ 0 \end{bmatrix}}
\end{equation}
\added{for the $i^{\rm th}$ neuron, where $I_{b,0}^i$ is the nominal constant baseline current.}

\begin{figure}[htb]
\begin{center}
\includegraphics[height=2.4in]{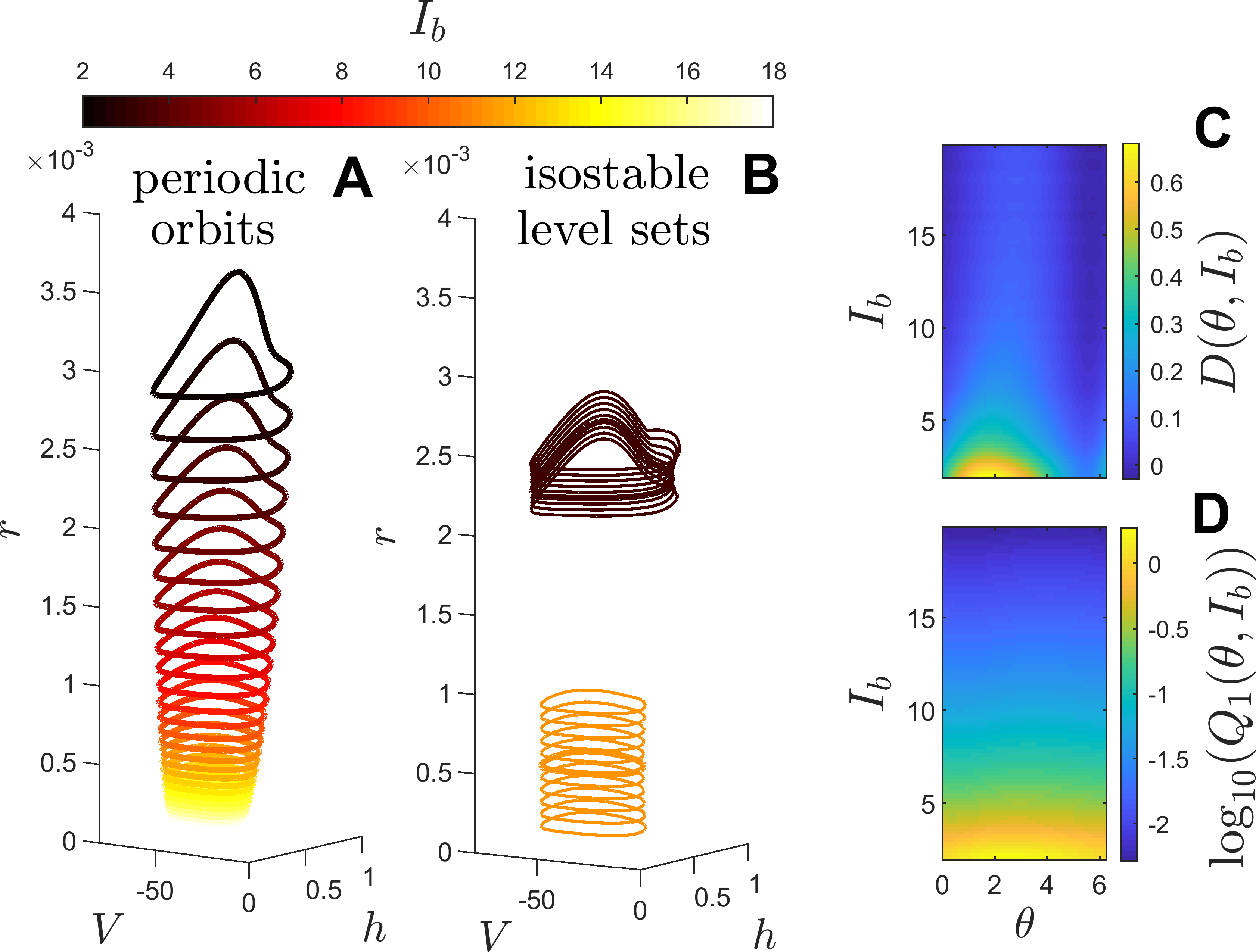}
\end{center}
\caption{Panel A gives examples of the periodic orbits \added{of \eqref{neurmod}} at various values of $I_b$ and Panel B  gives level sets of $\psi_1$ for two different $I_b$ values with the thick lines denoting the periodic orbits.  Panels C and D show $D(\theta,I_b)$ and $Q_1(\theta,p)$ that are used in \eqref{truncmain} to characterize how changes in $I_b$ influence the phase and isostable coordinates, respectively.}
\label{backgroundinfo}
\end{figure}

For the range of allowable $I_b$, the \added{principal} Floquet exponents are $\kappa_1 \in [-0.029,-0.019]$.  \added{The other Floquet exponents take values $\kappa_2 \in [-.45,-.29]$ and $\kappa_3 \in [-2.24,-2.09]$; these Floquet exponents are much larger in magnitude and will be truncated from the adaptive reduction \eqref{truncmain}.}  \added{In this example, $\theta(x,p) = 0$ corresponds to the moment that the $V$ reaches its peak value during an action potential on each of the limit cycles parameterized by the choice of $I_b$}.   \added{Compared to the full 4-dimensional model equations \eqref{neurmod}}, the resulting reduction has 3 variables for each neuron $\theta^i$, $\psi_1^i$, and $I_b^i$ with $i$ denoting the neuron number.    Because there is only one non-constant parameter and one isostable coordinate per neuron, a relatively simple choice of $G_p$ as suggested by Equation \eqref{firstgp} can be used:~$\dot{I}^i_b = G_p(I_b^i,\theta^i,\psi_1^i,) = - \alpha \psi_1^i \frac{\partial \psi_1^i}{\partial I^i_b},$ where $\alpha >0$ and $\partial \psi_1^i /\partial I_b^i = Q_1(\theta^i, I_b^i)+\psi_1 R_1^1(\theta^i, I_b^i)$.     Here $\alpha = 1000$ is chosen, although the exact choice is not important provided it keeps $\psi_1$ small.  General information for each reduced neuron from \eqref{neurmod} is shown in Figure \ref{backgroundinfo}.  

\begin{figure}[htb]
\begin{center}
\includegraphics[height=3in]{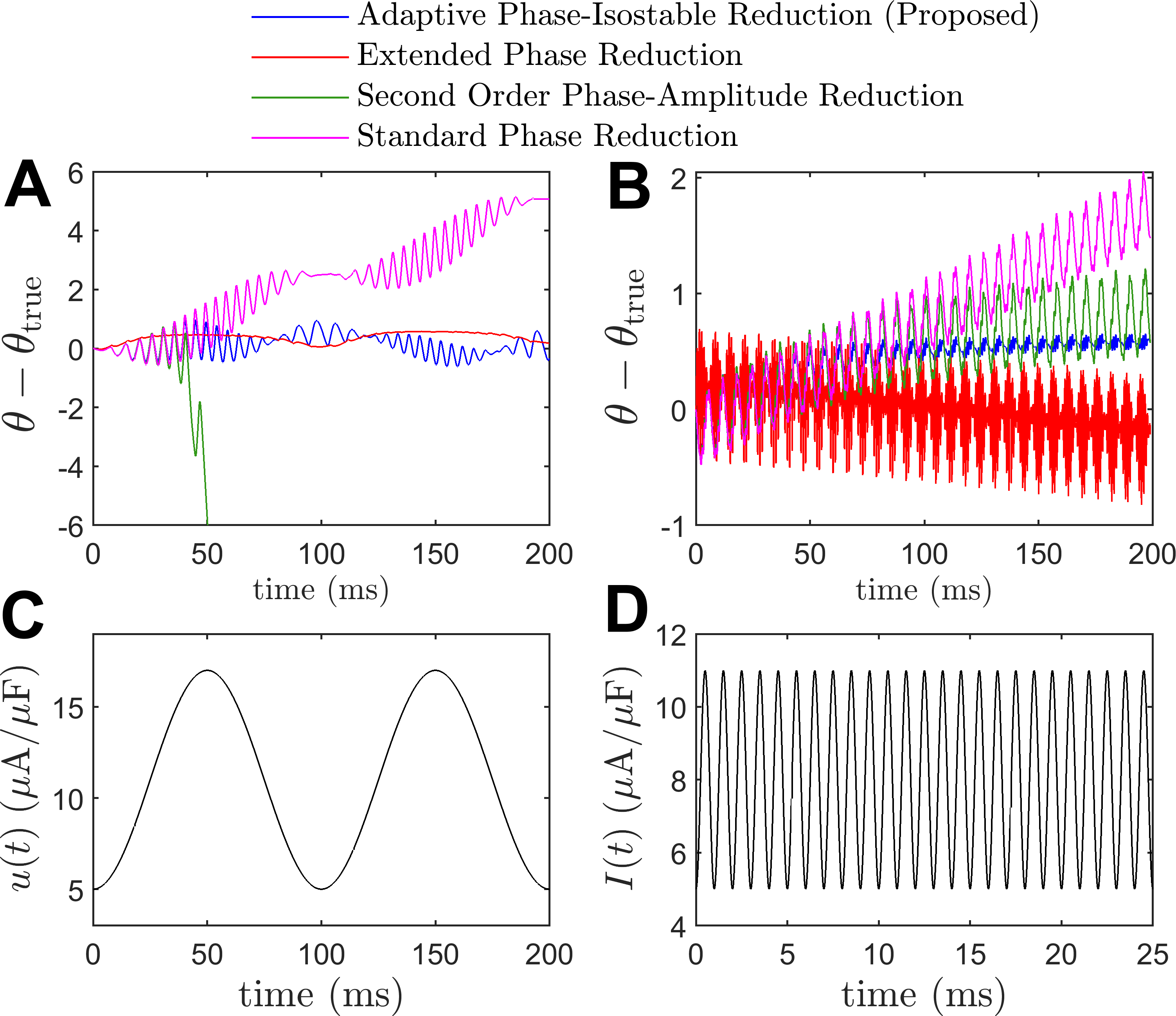}
\end{center}
\caption{For the slowly (resp.,~rapidly) varying input shown in panel \added{C} (resp.,~D), Figure A (resp.,~B) shows the error between the phase from a reduced model output and the phase of the full model equations.   \added{Here,  the proposed adaptive reduction strategy is considered and the performance is compared to three other reduction strategies.}    All initial conditions correspond to $\theta = 0$ on the limit cycle resulting from taking $I_b = 5 \mu {\rm A}/\mu {\rm F}$.}
\label{openloopresults}
\end{figure}

For the moment, a single neuron will be simulated in the absence of synaptic current \added{($g_{\rm syn} =0$)} and the proposed adaptive phase-isostable reduction strategy will be compared against three other reduction strategies using open loop sinusoidal input for $u(t)$.  The first is standard phase reduction \cite{erme10}, \cite{izhi07}, which takes the form $\dot{\theta}^i = \omega(I_{b,0}) + Z_V(\theta^i,I_{b,0}) (u(t)  - I_{b,0})$ where $I_{b,0} = 5 \mu {\rm A}/\mu {\rm F}$  is a constant baseline current and $Z_V$ is the component of $Z(\theta,I)$ in the $V$-direction.  \added{The second strategy is the second order accurate phase-amplitude reduction from \eqref{firstord} that is also implemented with respect to a constant baseline current input $I_{b,0} = 5 \mu {\rm A}/\mu {\rm F}$}.  The third strategy is the extended phase reduction proposed in \cite{kure13} and described by \eqref{extpred}.  This reduction assumes that the input is sufficiently slowly varying which yields a single equation for the phase dynamics $\dot{\theta} = \omega(u(t)) + D(\theta,u(t)) \cdot \dot{u}$.  Results are shown for both a slowly and more rapidly varying input in Figure \ref{openloopresults}.  The true phase, $\theta_{\rm true}$, of the full model simulations is approximated by taking $\theta = 0$ to correspond to the moment the transmembrane voltage reaches a maximum value during each cycle and is linearly interpolated at all points in-between.  The proposed adaptive phase-amplitude reduced model accurately characterizes the system behavior for both inputs unlike the other strategies.

\begin{figure}[htb]
\begin{center}
\includegraphics[height=2.7in]{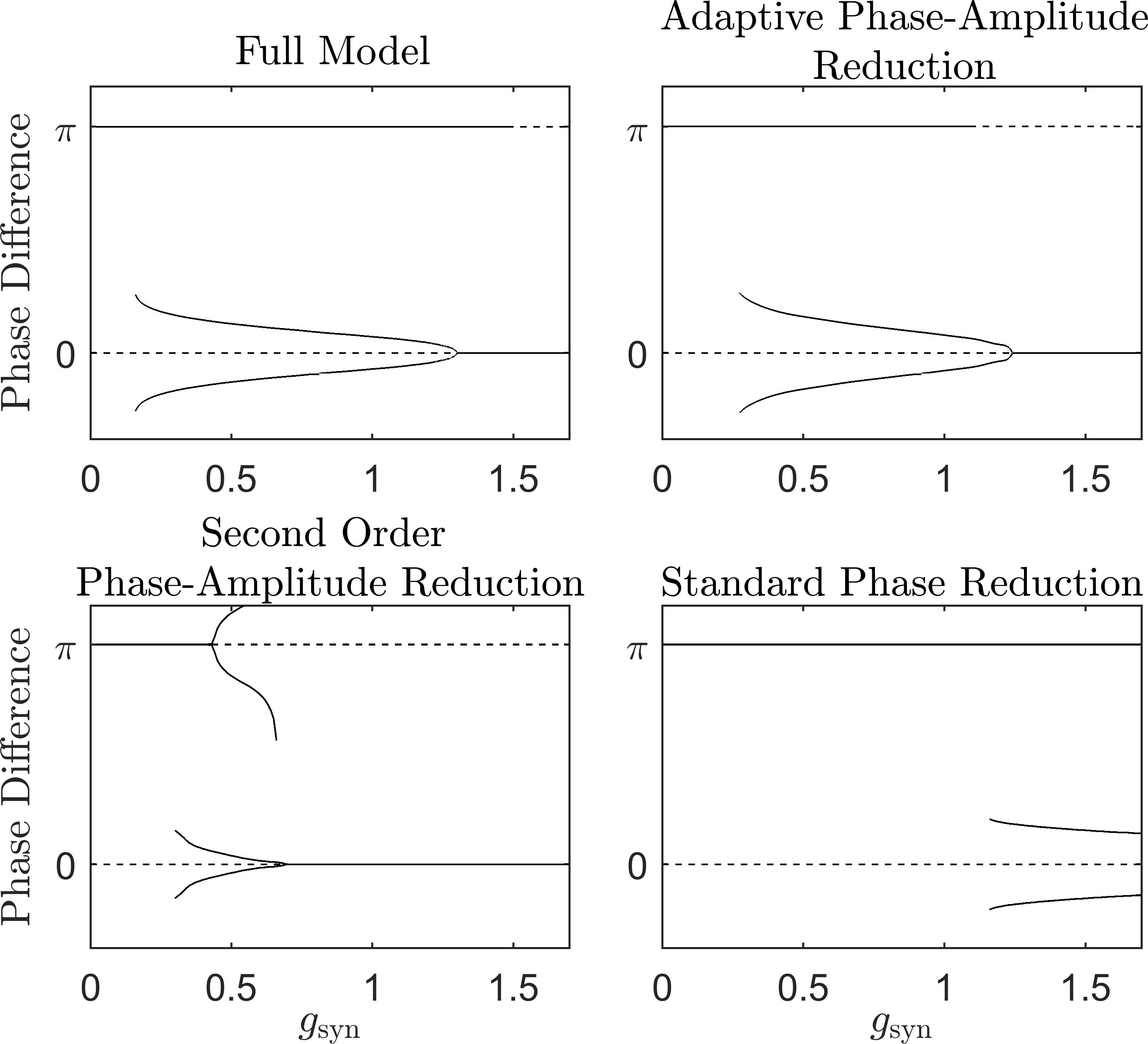}
\end{center}
\caption{Plots of the phase differences for stable (solid lines) and unstable (dashed lines) phase locked solutions for full simulations of \eqref{neurmod} with $n=2$ and three different reduction strategies.  The proposed adaptive phase-amplitude reduction accurately reflects bifurcations occurring in the full model for the full range of $g_{\rm syn}$ considered while the other two strategies are only accurate for small coupling strengths.}
\label{bifndiag}
\end{figure}

Next, synchronization of $N_{\rm neur} = 2$ identical, synaptically coupled neurons with nominal baseline current $I_{b,0} = 3.5 \mu {\rm A}/\mu {\rm F}$ is considered.  For these simulations $u(t)=0$ \added{(Note here that $U^i(t) \neq 0$ since the neurons are synaptically coupled)}.  Results from the full model \eqref{neurmod} are compared to the proposed adaptive phase-amplitude reduction \eqref{truncmain}, second order accurate phase-amplitude reduction from \eqref{firstord}, and standard phase reduction \eqref{predapx} \added{and shown in Figure \ref{bifndiag}}.    To implement the adaptive and second order phase-amplitude reduction, the variables $V_i$ and $w_i$ are approximated using the relationship \eqref{firstorderdx}.  These states are then used to compute the synaptic current $I_{\rm syn}(V_i,w_j)$.  The same procedure is used for the standard phase reduction with each isostable coordinate taken to be zero.   The adaptive phase-amplitude strategy accurately reflects the model behavior for large coupling strengths while the other two methods fail.  The extended phase reduced model from \cite{kure13} yields an unstable model when using this strategy to compute the synaptic current and cannot be considered.   \added{It is noted that the rise time of the synaptic variable is relatively slow and the coupling is excitatory so that that in the weak coupling limit, the existence of a stable antiphase state is in agreement with results presented by \cite{van94}.}


\FloatBarrier

\subsection{Entrainment for Populations of Coupled Oscillators} \label{circsec}
Here, a much higher dimensional model will be considered that describes the population oscillation of $N = 10$ coupled oscillators \cite{gonz05} which has been used to model populations of suprachiasmatic nucleus cells and resulting circadian rhythms
\begin{align} \label{circmodel}
\dot{a}_i &= h_1 \frac{K_1^n}{K_1^n + c_i^n}  -  h_2  \frac{a_i}{K_2+a_i} + h_c \frac{KF(t)}{K_c + KF(t)} + S_i L(t) + L_0,   \nonumber \\
\dot{b}_i &= h_3 a_i - h_4 \frac{b_i}{K_4+b_i},  \nonumber \\
\dot{c}_i &=  h_5 b_i - h_6\frac{c_i}{K_6+c_i} ,   \nonumber \\
\dot{d}_i &=     h_7 a_i - h_8 \frac{d_i}{K_8 + d_i}, \nonumber \\
\quad i &= 1,\dots,N.
\end{align}
Here variables $a_i$, $b_i$, and $c_i$ represent the concentrations of mRNA clock gene, associated protein, and nuclear form of the protein, respectively, for cell $i$, $d_i$ is a neurotransmitter that determines the mean-field coupling $F(t) \equiv (1/N)\sum_{j = 1}^N d_i(t)$.  $S_i = 1 + (i-1)/45$ denotes the sensitivity to light.  All other cell parameters are identical and taken to be the same as those from Figure 1 of \cite{gonz05} except for $n=7$, $h_1 = 1.05$, $h_2 = 0.525$ and $h_c = 0.2$.  For this choice of parameters, the nominal period of the population oscillation is 24.7 hours.    Two sigmoidal curves define a $T$-periodic light-dark cycle $L(t) =  L_{\rm max} /[{1+\exp(-4 (t_s -T/4))}]  -  L_{\rm max}/[{1+\exp(-4 (t_s -3T/4))} ] $ where $t_s = {\rm mod}(t,T)$ and $L_{\rm max}$ is the maximum light intensity. $L_0 \in [-0.010, 0.019]$ is treated as an \added{adaptive} parameter in the adaptive reduction.

\begin{figure}[htb]
\begin{center}
\includegraphics[height=4.0in]{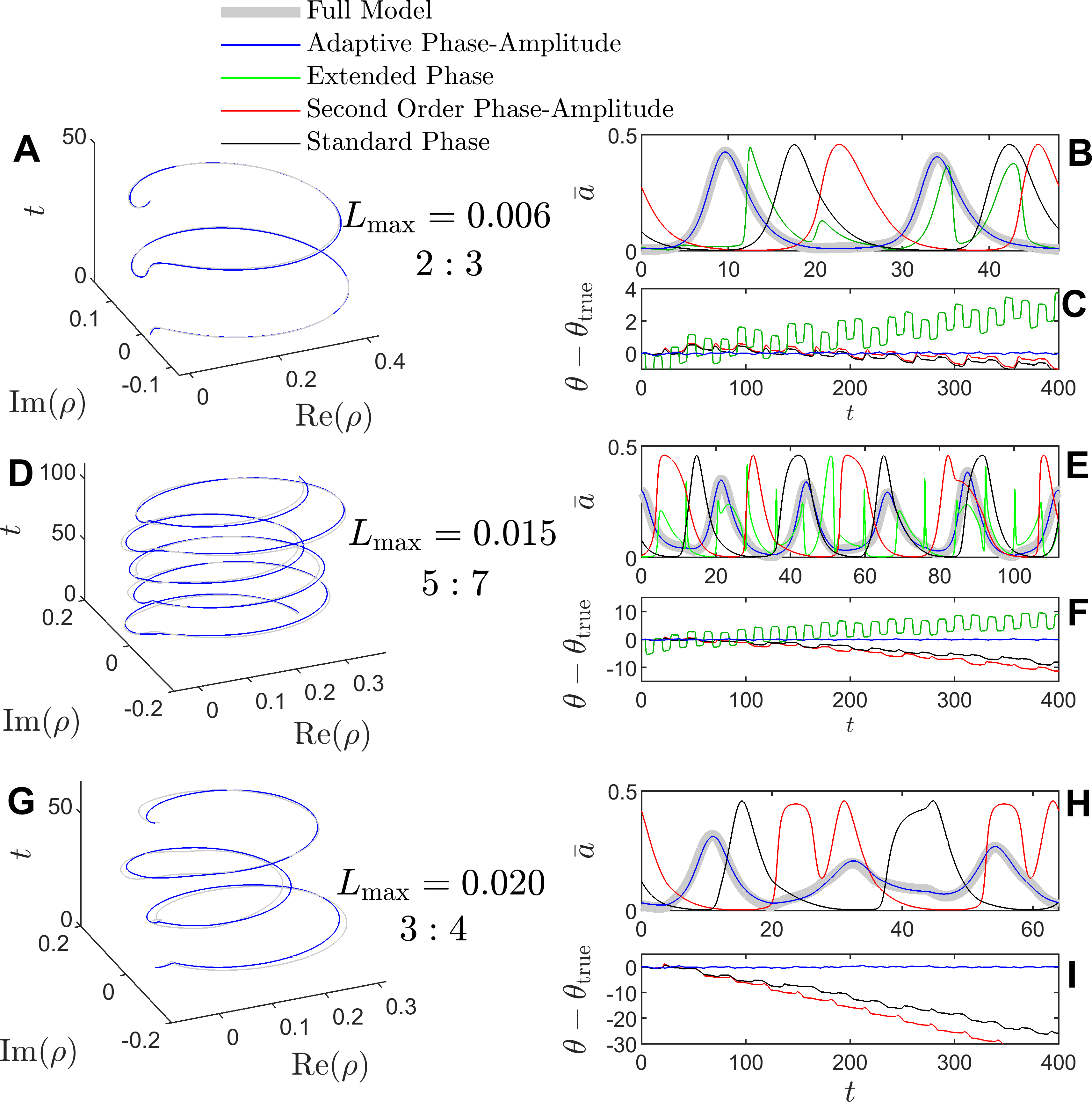}
\end{center}
\caption{Phase locking results with $N = 10$ oscillators in \eqref{circmodel} using various reduction strategies are illustrated when taking the period of $L(t)$ to be $T = 16$ hours.  The maximum light intensity $L_{\rm max}$ = 0.006, 0.015, and 0.020 results in 2:3, 5:7, and 3:4 locking, respectively, once steady state is achieved for the full model \eqref{circmodel}.  Corresponding results are shown in panels A-C, D-F, and G-I, respectively.  Panels B, E, and H show $\bar{a}(t) = (1/N) \sum_{i = 1}^N a_i(t)$ for the phase locked solutions of the full model and adaptive phase-amplitude reduced model.  None of the other reduced models exhibit phase locking for these choices of parameters, and simulated values are shown for reference.   The output is estimated using \eqref{firstorderdx}.  In reduction strategies that do not use isostable coordinates, each $\psi_j$ is taken to be zero when estimating the output from \eqref{firstorderdx}.   Panels A, D, and G show the real and imaginary components of $\rho(t) = \bar{a}(t) \exp(i \theta)$ for the phase locked solutions \added{plotted against time}.  The phase of the full model, $\theta_{\rm true}$ is approximated by taking $\theta_{\rm true} = 0$ each time a local maximum of $\bar{a}$ is achieved and interpolating the phase linearly for all times in-between.   Panels C, F, and I show the cumulative phase differences between the various phase reduced models and the full model over the course of each simulation.}
\label{circresult}
\end{figure}

For all allowable $L_0$, when $L(t) = 0$, \eqref{circmodel} has a stable synchronized orbit that can be studied in terms of a population oscillation.  Phase reduction has been applied for such systems previously, \cite{wils19phase}, \cite{kawa08}, \cite{leva10}, but the resulting models are typically only valid for prohibitively small inputs.  As shown here, the adaptive phase-amplitude reduction performs well compared to other methods when using particularly large inputs.   As explained in \cite{wils19phase}, because the population oscillation of \eqref{circmodel} results from mean-field coupling and because the inputs considered are of low rank, repeated isostable coordinates that arise due to symmetries can be ignored to first order accuracy in the isostable coordinates.  As a result, for a given $L_0$ and any $N \geq 2$,  the population oscillation of \eqref{circmodel} can be characterized with a reduced order model of the form \eqref{firstordadapt} that only requires 3 isostable coordinates and 1 phase coordinate.  \added{A reduction of a similar model was considered in \cite{wils19phase} where more details about this preliminary reduction and the specific choice isostable coordinates can be found}.   Including the time-varying parameter $L_0$, the adaptive phase-amplitude reduction \eqref{truncmain} requires a total of 5 variables compared to $4N$ in the full model. \added{In this example, $\theta(x,p) = 0$ corresponds to the moment that the $a_1$ reaches its peak value on each of the limit cycles parameterized by the choice of $L_0$.}   Additionally, the equation describing the adaptive variable $L_0$ is $\dot{L}_0 = -\alpha  \sum_{i = 1}^3  (  \psi_i  \frac{\partial \psi_i}{\partial L_0})$  with $\alpha = 0.0001$.   \added{This choice of $G_p$ is inspired by the formulation suggested by Equation \eqref{firstgp}} and is sufficient to keep the magnitude of the isostable coordinates small in the results to follow.     \added{In order to provide a direct comparison with \cite{wils19phase},} the $C_j^k(\theta,p_0)$ terms are ignored in the reduction.  \added{Consequently, the isostable coordinate dynamics are accurate to first order in the basis of the non-truncated isostable coordinates.}  Figure \ref{circresult} shows phase locking results for a population of $N = 10$ oscillators when the period of external light forcing is $T = 16$ hours.   Each initial condition in simulations of \eqref{circmodel} starts from phase $\theta = 0$ on the periodic orbit (i.e.,~taking all isostable coordinates equal to zero).   Results are shown in Figure \ref{circresult} for multiple types of phase reductions.  The standard phase reduction \eqref{predapx} and the second order phase reduction of the form \eqref{firstord} are implemented using $L_0 = 0$.  The extended phase reduction from \cite{kure13} is implemented by taking $L_0 = 0$ and viewing $L(t)$ as a slowly varying parameter yielding an equation of the form $\dot{\theta} = \omega(L) + D(\theta,L) \dot{L}$ (see also Equation \eqref{extpred} from Section \ref{apxb}).   Once again, the proposed adaptive phase-amplitude reduction is the only reduction framework that accurately reflects the model behavior for the inputs considered.

\section{Conclusion} \label{concsec}
While many techniques exist for analyzing oscillatory dynamical systems in the weakly perturbed limit, existing techniques have struggled to characterize behavior in strongly perturbed systems.  The method proposed in this work provides a general, adaptive phase-amplitude reduction framework for strongly perturbed oscillatory dynamical systems.   In contrast to other phase and phase-amplitude reduction techniques, this method works for arbitrary inputs; no $\mathcal{O}(\epsilon)$ constraints are placed on either the magnitude or rate of change of allowable inputs.   In the numerical examples presented here, the proposed framework accurately reflects the behavior of the underlying models in regimes where all other considered phase and phase-amplitude reductions fail.  

\added{It should be emphasized that many prior phase and phase-amplitude reduction frameworks (such as those described in Section \ref{apxb}) generally yield lower dimensional models than the proposed adaptive reduction framework; they also generally and require less computational effort to identify the necessary terms of the reductions.  As such, the proposed adaptive reduction framework method will likely find use in situations where these other reduction frameworks do not accurately replicate the full order model dynamics.}


\added{To implement the proposed adaptive reduction strategy, a family of limit cycles is considered that emerge as an adaptive parameter set is changed.  This adaptive parameter set is continuously adjusted so that the underlying system state stays close at all times to the limit cycle associated with current choice of the  adaptive parameters.  This strategy can be viewed as an extension of the approach described in Equation \eqref{extpred} from \cite{kure13} (cf.,~\cite{park16}) that also considers a similar family of limit cycles.  However, because the strategy proposed in \cite{kure13} does not consider the amplitude coordinates of each limit cycle, allowable inputs are explicitly required to vary sufficiently slowly.  Indeed, considering the results presented in Figure \eqref{openloopresults}, the reduction approach from Equation \eqref{extpred} works well when a low frequency input is applied but fails for a high frequency input.  By actively considering the amplitude coordinates and adjusting the nominal parameter set as appropriate, the proposed adaptive reduction strategy of the form \eqref{truncmain} accurately captures system behavior regardless of the forcing frequency.}


\added{This work provides a foundation for the implementation of an adaptive reduction framework, but there are still many questions and limitations remaining.  Foremost,  per Assumption C from Section \ref{assumpsec}, in order to implement the adaptive reduction one must identify a parameter update rule $\dot{p} = G_p(p,\theta,\psi_1,\dots,\psi_\beta)$ that ensures that the isostable coordinates remain small at all times.  Section \ref{designgp} discusses some general design heuristics that can be used in situations where there is only one adaptive parameter.  However, situations where more than one adaptive parameter is required are not considered in this manuscript.  The design of an appropriate $G_p$ could be posed as a nonlinear control problem and will be the subject future investigations.   Furthermore, the proposed method requires a sufficient degree of smoothness in the underlying system equations as mandated by assumptions A and B from Section \ref{assumpsec}.  These assumptions exclude, for instance, the use of the adaptive reduction framework in systems with discontinuities such as the piecewise smooth systems considered in \cite{park16ar}.  It is likely that the proposed adaptive reduction framework could be extended for use in these piecewise smooth systems systems but such applications are not considered here.   Finally, it may also be possible to extend the proposed method to include critical points of bifurcations in the allowable parameter set (e.g.,~transitioning from stable limit cycle to a stable fixed point solutions through a Hopf bifurcation).  }

This work was supported by National Science Foundation Grant No.~CMMI-1933583.

\begin{appendices}
\renewcommand{\thetable}{A\arabic{table}}  
\renewcommand{\thefigure}{A\arabic{figure}} 
\renewcommand{\theequation}{A\arabic{equation}} 
\setcounter{equation}{0}
\setcounter{figure}{0}

\section{Model Reduction Using Adaptive Phase-Amplitude Coordinates}  \label{apxa}
Consider a general ordinary differential equation
\begin{equation} \label{underlyingmodel}
\dot{x} = F(x,p_0) + U(t),
\end{equation}
where $x \in \mathbb{R}^N$ denotes the state, $p_0\in \mathbb{R}^M$ is a constant collection of nominal parameters, $F$ gives the dynamical behavior, and $U(t)$ represents an external perturbation.   \added{To implement the proposed adaptive phase-amplitude reduction, consider the shadow system from \eqref{shadoweq}.  The adaptive parameters are chosen and an associated function $G_p$ is defined so that Assumptions A-E from Section \ref{assumpsec} are satisfied}.  If no isostable coordinates are truncated in the resulting adaptive phase-amplitude reduction from \eqref{truncmain}, the equations are
  \begin{align} \label{extendedfirst}
\dot{\theta} &= \omega(p) + \left(Z(\theta,p)  + \sum_{k = 1}^{N-1} \psi_k B^k(\theta,p) \right)U_e(t,p,x)  +  \left( D(\theta,p) + \sum_{k=1}^{N-1} \psi_k E^k(\theta,p)  \right) \dot{p} \nonumber \\
&   + \mathcal{O}(\psi_1^2) + \dots + \mathcal{O}(\psi_{N-1}^2), \nonumber \\
\dot{\psi}_j &= \kappa_j(p) \psi_j +  \left(I_j(\theta,p)  + \sum_{k = 1}^{N-1} \psi_k C_j^k(\theta,p) \right) U_e(t,p,x) +  \left( Q_j(\theta,p) + \sum_{k=1}^{N-1} \psi_k R_j^k(\theta,p)  \right) \dot{p}, \nonumber \\
&  + \mathcal{O}(\psi_1^2) + \dots + \mathcal{O}(\psi_{N-1}^2), \nonumber \\
& j = 1, \dots, N-1, \nonumber \\
\dot{p} &= G_p(p,\theta,\psi_1,\dots,\psi_{N-1}).
\end{align}
Here, the  $\mathcal{O}(\psi_1^2) + \dots + \mathcal{O}(\psi_{N-1}^2)$ terms result from truncation of the gradients of the phase and isostable coordinates after first order accuracy.  \added{Recall that $U_e$ is also a function of $p_0$, but this dependence is suppressed in Equation \eqref{extendedfirst} for clarity for exposition.} Supposing some $G_p$ could be obtained that could keep all of the isostable coordinates small, \eqref{extendedfirst} would be accurate for describing the phase and amplitude dynamics, however, \eqref{extendedfirst} would have $N+M$ states, which is more than the original system started with.  Here, it will be shown that provided some of the Floquet exponents from \eqref{underlyingmodel} are sufficiently large in magnitude, the associated isostable coordinates can be ignored.

\added{Drawing on Assumption D from Section \ref{assumpsec},} the Floquet exponents $\kappa_{\beta+1}(p), \dots, \kappa_{N-1}(p)$  are $\mathcal{O}(1/\epsilon)$ terms for all values of $p$ where $0 \leq \epsilon \ll 1$, i.e.,~their associated isostable coordinates rapidly decay in response to perturbations.  \added{As such, the rapidly and slowly decaying isostable coordinates will be considered separately rewriting Equation \eqref{extendedfirst} as}
\begin{align} \label{sepeq}
\dot{\theta} &= \omega(p) + \left(Z(\theta,p)  + \sum_{k = 1}^{\beta} \psi_k B^k(\theta,p)  +  \sum_{k = \beta+1}^{N-1} \psi_k B^k(\theta,p)  \right) \cdot U_e(t,p,x)  \nonumber \\
& \quad + \left( D(\theta,p) + \sum_{k=1}^{\beta} \psi_k E^k(\theta,p)  +  \sum_{k=\beta+1}^{N-1} \psi_k E^k(\theta,p)  \right) \cdot \dot{p},   \nonumber \\
\dot{\psi}_j &= \kappa_j(p) \psi_j +  \left(I_j(\theta,p)  + \sum_{k = 1}^{\beta} \psi_k C_j^k(\theta,p)    +  \sum_{k = \beta+1}^{N-1} \psi_k C_j^k(\theta,p)   \right) \cdot U_e(t,p,x)  \nonumber \\
& \quad  +  \left( Q_j(\theta,p) + \sum_{k=1}^{\beta} \psi_k R_j^k(\theta,p)  + \sum_{k=\beta+1}^{N-1} \psi_k R_j^k(\theta,p)    \right) \cdot \dot{p}, \nonumber \\
& j = 1, \dots, \beta, \nonumber \\
\dot{\psi}_i &= \kappa_i(p) \psi_i +  \left(I_i(\theta,p)  + \sum_{k = 1}^{\beta} \psi_k C_i^k(\theta,p)   +  \sum_{k = \beta+1}^{N-1} \psi_k C_i^k(\theta,p)   \right) \cdot U_e(t,p,x)  \nonumber \\
& \quad  +  \left( Q_i(\theta,p) + \sum_{k=1}^{\beta} \psi_k R_i^k(\theta,p)  + \sum_{k=\beta+1}^{N-1} \psi_k R_i^k(\theta,p)    \right) \cdot \dot{p}, \nonumber \\
& i = \beta+1, \dots, N-1, \nonumber \\
\dot{p} &= G_p(p,\theta,\psi_1,\dots, \psi_\beta),
\end{align}
where  \added{$G_p$ no longer depends on the rapidly decaying isostable coordinates}.  Above, the order of accuracy of the neglected terms has been omitted for clarity of exposition.  Above, the rapidly decaying isostable coordinates have been explicitly separated from the other isostable coordinates.  Also,  $\dot{p}$ is no longer a function of the rapidly decaying isostable coordinates  in \eqref{sepeq}.  Suppose that at an initial time $t = 0$, each $\psi_j = 0$, i.e., the solution is on the periodic orbit $x^\gamma_p$.   \added{Also recall from Assumption E in Section \ref{assumpsec}} that the norms of $U_e(t,p,x)$ and $G_p(p,\theta,\psi_1,\dots,\psi_\beta)$ (i.e.,~that determines $\dot{p})$ are bounded for all time and $p$ so that $||U_e(t,p,x)||_1 \leq M_U$ and  $||  \dot{p} ||_1 \leq M_p$ where $|| \cdot ||_1$ denotes the $1$-norm.  Taking the absolute value of any $\psi_i$ for $i \geq \beta+1$ from \eqref{sepeq}, \added{recalling that $\kappa_i(p)<0$}, one can show
\begin{align} \label{longleq}
| \dot{\psi}_i| &\leq  \kappa_i(p) |\psi_i| + \max_{\theta,p,t} |  I_i(\theta,p)\cdot U_e(t,p,x) | +  \beta  \max_{\substack{\theta,p,t\\k = 1,\dots,\beta}}| C_i^k (\theta,p) \cdot U_e(t,p,x)|  \max_{k = 1,\dots,\beta}  |\psi_k| \nonumber \\
  & \quad +    (N-\beta-1)  \max_{\substack{  \theta,p,t \\ k = \beta+1, \dots,N-1  }}| C_i^k (\theta,p) \cdot U_e(t,p,x)  | \max_{k = \beta+1,\dots,N-1} |  \psi_k  | \nonumber \\
 & \quad + \max_{\theta, p, \dot{p}} |  Q_i(\theta,p) \cdot \dot{p}  | +  \beta \max_{\substack{  \theta,p, \dot{p}\\ k = 1, \dots, \beta  }} |   R_i^k(\theta,p) \cdot \dot{p}|   \nonumber \\
 &  \quad +    (N-\beta-1)  \max_{\substack{  \theta,p, \dot{p} \\ k = \beta+1, \dots,N-1  }}  |   R_i^k(\theta,p) \cdot \dot{p}   |  \max_{k = \beta+1,\dots,N-1} |  \psi_k  |  \nonumber \\
 &\leq  \max_p (\kappa_i(p))  | \psi_i | + M_U \max_{\theta,p} ||  I_i(\theta,p) ||_1 + \beta M_U \max_{\substack{\theta,p\\k = 1,\dots,\beta}}||  C_i^k(\theta,p)||_1 \max_{k = 1,\dots,\beta}  |\psi_k| \nonumber \\
 & \quad + (N-\beta-1) M_U   \max_{\substack{  \theta,p \\ k = \beta+1, \dots,N-1  }}  ||C_i^k(\theta,p) ||_1  \max_{k = \beta+1,\dots,N-1} |  \psi_k  |   \nonumber \\
 & \quad + M_p \max_{\theta,p} || Q_i (\theta,p)|| + \beta M_p \max_{\substack{  \theta,p\\ k = 1, \dots, \beta  }} | |  R_i^k(\theta,p)||_1     \max_{k = 1,\dots,\beta}  |\psi_k| \nonumber \\
 & \quad +   (N-\beta-1) M_p  \max_{\substack{  \theta,p \\ k = \beta+1, \dots,N-1  }}  ||   R_i^k(\theta,p) ||_1  \max_{k = \beta+1,\dots,N-1} |  \psi_k  | .
\end{align}
\textcolor{black}{Continuing to simplify \eqref{longleq} and rearranging terms yields 
\begin{align} \label{simpterms}
| \dot{\psi}_i| &\leq  \max_p (\kappa_i(p))   |\psi_i| + M_U a_1 + M_U a_2 \max_{k = 1,\dots,\beta}  |\psi_k| + M_U a_3    \max_{k = \beta+1,\dots,N-1} |  \psi_k  |   \nonumber \\
 & \quad + M_p a_4 +  M_p   a_5  \max_{k = 1,\dots,\beta}  |\psi_k| +    M_p  a_6  \max_{k = \beta+1,\dots,N-1} |  \psi_k  |  ,
\end{align}}
where $a_1 \equiv \max_{\theta,p,i = \beta+1,\dots,N-1} ||  I_i(\theta,p) ||_1$, $a_2 \equiv \beta \max_{\substack{\theta,p,k = 1,\dots,\beta},i = \beta+1,\dots,N-1}||  C_i^k(\theta,p)||_1$, $a_3 \equiv (N-\beta-1)   \max_{\substack{  \theta,p, k = \beta+1, \dots,N-1,i = \beta+1,\dots,N-1  }}  ||C_i^k(\theta,p) ||_1  $,  $a_4 \equiv \max_{\theta,p,i = \beta+1,\dots,N-1} || Q_i (\theta,p)||$, $a_5 \equiv  \beta \max_{\substack{  \theta,p, k = 1, \dots, \beta,i = \beta+1,\dots,N-1  }} | |  R_i^k(\theta,p)||_1 $, and $a_6 \equiv (N-\beta-1)  \max_{\substack{  \theta,p , k = \beta+1, \dots,N-1 ,i = \beta+1,\dots,N-1 }}  ||   R_i^k(\theta,p) ||_1$.  \added{Note that $a_1$ through $a_6$ are related to the dimension of the system and the terms of the partial derivatives.  Here, the magnitude of the neglected Floquet exponents must be large enough so that $\epsilon$ is small enough for $a_1$ through $a_6$ to be considered $\mathcal{O}(1)$ terms and not $\mathcal{O}(1/\epsilon)$ terms.}  \textcolor{black}{Recalling that $ \max_p (\kappa_i(p))   < 0$ and considering \eqref{simpterms}, notice that if 
\begin{align}
|\psi_i| > \frac{1}{- \max_p (\kappa_i(p))  } \bigg[& M_U a_1 + M_U a_2 \max_{k = 1,\dots,\beta}  |\psi_k| + M_U a_3    \max_{k = \beta+1,\dots,N-1} |  \psi_k  |   \nonumber \\
 & \quad + M_p a_4 +  M_p   a_5  \max_{k = 1,\dots,\beta}  |\psi_k| +    M_p  a_6  \max_{k = \beta+1,\dots,N-1} |  \psi_k  |  \bigg],
\end{align}
is true, then consequently $| \dot{\psi}_i| < 0$.  Accordingly, 
\begin{align}\label{upperbound}
 - \max_p (\kappa_i(p))   |\psi_i| &\leq M_U a_1 + M_U a_2 \max_{k = 1,\dots,\beta}  |\psi_k| + M_U a_3    \max_{k = \beta+1,\dots,N-1} |  \psi_k  |   \nonumber \\
 & \quad + M_p a_4 +  M_p   a_5  \max_{k = 1,\dots,\beta}  |\psi_k| +    M_p  a_6  \max_{k = \beta+1,\dots,N-1} |  \psi_k  |  
\end{align}
gives an upper bound for each $|\psi_i|$.  }  Considering that \eqref{upperbound} is valid for any isostable coordinate $\psi_i$ with $i \geq \beta+1$, particularly for the $\psi_i$ with the maximum absolute value, one can continue to manipulate the bounds on the left side of \eqref{upperbound} to yield
\begin{align}
\min_{\substack{p\\i = \beta+1,\dots,N-1}}(-\kappa_i(p))    \max_{k = \beta+1,\dots,N-1} |  \psi_k  | & \leq   - \max_p (\kappa_i(p))    \max_{k = \beta+1,\dots,N-1} |  \psi_k  |  \nonumber \\
& \leq  M_U a_1 + M_U a_2 \max_{k = 1,\dots,\beta}  |\psi_k| + M_U a_3    \max_{k = \beta+1,\dots,N-1} |  \psi_k  |   \nonumber \\
 & \quad + M_p a_4 +  M_p   a_5  \max_{k = 1,\dots,\beta}  |\psi_k| +    M_p  a_6  \max_{k = \beta+1,\dots,N-1} |  \psi_k  |  ,
\end{align}

Recalling that each $\kappa_i$ is an order $1/\epsilon$ term for $i = \beta + 1, \dots, N-1$, $\min_{\substack{p,i = \beta+1,\dots,N-1}}(-\kappa_i(p)) -  M_U a_3 - M_p a_6 > 0$, so that
\begin{equation} \label{psibound}
\max_{k = \beta+1,\dots,N-1} |  \psi_k  |  \leq   \frac{M_U a_1 + M_pa_4 +  (M_U a_2 + M_p a_5) \max_{k = 1,\dots,\beta}  |\psi_k|  }{\min_{\substack{p,i = \beta+1,\dots,N-1}}(-\kappa_i(p)) -  M_U a_3 - M_p a_6  }.
\end{equation} 
Recalling that $\min_{\substack{p,i = \beta+1,\dots,N-1}}(-\kappa_i(p))$ is assumed to  be  $\mathcal{O}(1/\epsilon)$, Equation \eqref{psibound} implies that each $\psi_k$ must be an $\mathcal{O}(\epsilon)$ term.

With an explicit bound on the rapidly decaying isostable coordinates, it is possible to ignore them from Equation \eqref{sepeq}.  Doing so yields 
\begin{align} \label{trunceq}
\dot{\theta} &= \omega(p) + \left(Z(\theta,p)  + \sum_{k = 1}^{\beta} \psi_k B^k(\theta,p)    \right) \cdot U_e(t,p,x)   + \left( D(\theta,p) + \sum_{k=1}^{\beta} \psi_k E^k(\theta,p)    \right) \cdot \dot{p}   \nonumber \\
&+ \mathcal{O}(\epsilon) + \mathcal{O}(\psi_1^2) + \dots + \mathcal{O}(\psi_\beta^2), \nonumber \\
\dot{\psi}_j &= \kappa_j(p) \psi_j +  \left(I_j(\theta,p)  + \sum_{k = 1}^{\beta} \psi_k C_j^k(\theta,p)    \right) \cdot U_e(t,p,x)   +  \left( Q_j(\theta,p) + \sum_{k=1}^{\beta} \psi_k R_j^k(\theta,p)     \right) \cdot \dot{p} \nonumber \\
&+ \mathcal{O}(\epsilon) + \mathcal{O}(\psi_1^2) + \dots + \mathcal{O}(\psi_\beta^2), \nonumber \\
j &= 1, \dots, \beta, \nonumber \\
\dot{p} &= G_p(p,\theta,\psi_1,\dots, \psi_\beta),
\end{align}
where the $\mathcal{O}(\epsilon)$ terms come from the truncation of the rapidly decaying isostable coordinates.  Finally,  as long as a $G_p(p,\psi_1,\dots,\psi_\beta)$ can be designed so that each $\psi_1, \dots,\psi_\beta$ remains an $\mathcal{O}( \sqrt \epsilon )$ term for all  time $t >0$, the reduction \eqref{trunceq} is valid to leading order $\epsilon$.   As emphasized in the main text, Equation \eqref{trunceq} does not place any $\mathcal{O}(\epsilon)$ constraints on the input $U(t,p,x)$.

\end{appendices}


\end{document}